\documentclass[11pt]{article}
\usepackage[hypertexnames=false]{hyperref}
\usepackage{color}
\usepackage{bbm}
\usepackage{bm}
\usepackage[dvipsnames]{xcolor}
\usepackage{latexsym}
\usepackage{dsfont}
\usepackage{appendix}
\usepackage{amssymb}
\usepackage{subcaption}
\usepackage[margin=.9in]{geometry}
\usepackage{float}
\usepackage{algorithm,algorithmic}
\usepackage[numbers]{natbib}
\usepackage{amsmath, amsfonts,amssymb,euscript,array,enumerate,amsfonts,mathrsfs}
\usepackage[dvipsnames]{xcolor}
\usepackage{amsthm}
\newtheorem{Theorem}{Theorem}[section]
\newtheorem{Definition}[Theorem]{Definition}
\newtheorem{Proposition}[Theorem]{Proposition}
\newtheorem{Assumption}[Theorem]{Assumption}
\newtheorem{Lemma}[Theorem]{Lemma}
\newtheorem{Corollary}[Theorem]{Corollary}
\newtheorem{Remark}[Theorem]{Remark}
\newtheorem{Example}[Theorem]{Example}
\usepackage[T1]{fontenc}
\usepackage[latin1]{inputenc}
\usepackage{hhline}
\usepackage{graphicx}
\usepackage{verbatim}
\usepackage{eurosym}
\usepackage{dsfont}
\allowdisplaybreaks

\DeclareMathOperator{\E}{\mathbb{E}}

\DeclareMathOperator{\R}{\mathbb{R}}

\DeclareMathOperator*{\argmax}{arg\,max}

\makeatletter
\@addtoreset{equation}{section}

\newcommand*{\rom}[1]{\expandafter\@slowromancap\romannumeral #1@}
\makeatother
\newcommand{\vertiii}[1]{{\left\vert\kern-0.25ex\left\vert\kern-0.25ex\left\vert #1 
    \right\vert\kern-0.25ex\right\vert\kern-0.25ex\right\vert}}
\def\A{{\mathcal A}}
\def\S{{\mathcal S}}
\def\P{{\mathbb P}}
\def\PP{{\mathcal P}}
\def\E{{\mathbb E}}

\def\L{{\mathcal L}}
\def\l{{\ell}}

\def\M{{\mathcal M}}

\begin{document}
\title{Optimization frameworks and sensitivity analysis of Stackelberg mean-field games}

 \author{
Xin Guo\thanks{Department of Industrial Engineering and Operations Research, University of California, Berkeley, Berkeley, California, USA. Email: xinguo@berkeley.edu} 
\and Anran Hu\thanks{Mathematical Institute, University of Oxford, Oxford, United Kingdom. 
 Email: Anran.Hu@ maths.ox.ac.uk}
\and Jiacheng Zhang
 \thanks{Department of Industrial Engineering and Operations Research, University of California, Berkeley, Berkeley, California, USA
 Email:jiachengz@berkeley.edu }}
\date{}

\maketitle
\begin{abstract}
This paper proposes and studies a class of discrete-time finite-time-horizon Stackelberg mean-field games, with one leader and  an infinite number of identical and indistinguishable followers. In this game, the objective of the leader is to maximize her reward considering the worst-case cost over all possible $\epsilon$-Nash equilibria among followers.

A new analytical paradigm is established by showing the equivalence between this Stackelberg mean-field game  and a minimax optimization problem. 
 This optimization framework facilitates studying both analytically and numerically the set of Nash equilibria for the game; and leads to the  sensitivity and the robustness analysis of the game value. In particular, when there is model uncertainty, the game value for the leader suffers non-vanishing sub-optimality as the perturbed model converges to the true model. In order to obtain a near-optimal solution, the leader needs to be more pessimistic with anticipation of model errors and adopts a relaxed version of the original Stackelberg game.

\end{abstract}

\section{Introduction}\label{sec:intro}
Stackelberg game is a class of games with a bi-level structure, consisting of  two types of players: a leader and her followers. In this game, the leader takes an action first, then the followers  choose their strategies in response to the leader's action. The objective of the game  is for the leader to find the best action that maximizes her rewards, assuming that the followers always reach their best responses to her action. A special case of the Stackelberg game is the  principal-agent problem, where followers/agents maximize their rewards by following instructions proposed by the leader/principal; see \cite{holmstrom1987aggregation}, \cite{sannikov2008continuous}. 

Stackelberg game is not, however, restricted to only the principal-agent game setting. In the general case of a Stackelberg game, the leader and followers may have conflicting interests and followers may choose to deviate from leader's instruction. For instance, in an ad auction game, the interests of the ad platform and advertisers may indeed diverge: the ad platform takes actions  first by designing the mechanism of the auction game for the advertisers; the advertisers then make decisions to bid the ad slot to maximize their own profits. Unlike the  principal-agent problem where the principal can instruct the agents to a particular game equilibrium favorable for the principal \cite{elie2019tale},  the ad platform does not have control over the advertisers. 
Other instances of such Stackelberg games 
include resource allocation \cite{conitzer2002complexity},  optimal taxation \cite{zheng2020ai}, and 
security games \cite{korzhyk2011stackelberg},\cite{tambe2011security}.

Stackelberg games with such potentially deviating followers 
are not well understood and studied, despite its potentially wide range of applications. In addition, most theoretical studies of Stackelberg games to-date focus on the case with one leader and one follower. In this context,  \cite{conitzer2006computing} and \cite{letchford2010computing} study the computational complexities when the leader has full observation of the payoff function of  followers; and  \cite{fiez2019convergence} establishes the local convergence of first-order algorithms with an exact feedback. More recent works have also extended the study to the learning setting. In particular, 
\cite{blum2014learning},\cite{peng2019learning}, \cite{letchford2009learning} design learning algorithms assuming that the leader can induce followers' best responses under the best response oracle; \cite{bai2021sample} designs learning algorithms under a bandit setting and analyzes the corresponding sample complexities. A central tool used in the analysis of these learning algorithms is the stability/perturbation analysis of Stackelberg games. 
In contrast, the study of Stackelberg games with multiple followers is largely restricted to empirical studies and applications \cite{trejo2018adapting,jiang2020multi}, or the analysis of special settings such as stage games \cite{conitzer2006computing,coniglio2020computing}, extensive-form games \cite{letchford2010computing} and games with myopic players \cite{zhong2021can}. A recent work \cite{alcantara2020repeated} studies reinforcement learning of general Stackelberg games with multiple followers and establishes the convergence of the fitted model, yet without studying the performance of the learned policy due to the lack of stability analysis. This paper hence focuses on providing mathematical framework and answering the fundamental questions for Stackelberg games with more than one followers that remain open, 
including the regularity of the game value, and the stability of such games. 
These are the primary focus of this paper.

There are reasons why analyzing a Stackelberg game can be daunting.
First of all, there are typically a large number of followers, for instance in the ad auction where there could be tens of hundreds of advertisers. In order to study the followers' best response to the leader's action, it is necessary to analyze such  multiple-player games, which is notoriously  difficult. Secondly, it is challenging to analyze the best game strategy such as the Nash equilibrium when it is non-unique. In a Stackelberg game, it is essential to analyze leader's best action and its associated  value function based on the analysis of different  game equilibria of followers. Conceivably, analysis and results of Stackelberg games with one leader and one follower is of little help.  Meanwhile, there are almost no existing mathematical techniques to derive all Nash or game equilibria for non-zero-sum stochastic games. Current approaches rely on restricted assumptions such as contractivity (\cite{guo2019learning,guo2020general}) or monotonicity (\cite{perrin2020fictitious,perolat2021scaling}) even for mean-field-games, the simplified version of non-zero-sum stochastic games. The only exception is the recent work \cite{guo2022mf}, which establishes the equivalence between deriving an Nash equilibrium in an mean-field game and solving an optimization problem with smooth constraints, with the potential to analyze the set of all possible Nash equilibria.

\paragraph{Our work.}
In this paper, we consider a discrete-time finite-time-horizon Stackelberg mean-field games. In this game, there is a leader and there are an infinite number of identical and indistinguishable followers; the leader may not  have any control over the response from the followers once she takes an action, and her objective is to maximize her reward considering the worst-case cost over all possible equilibria among followers.

Given the possibility of multiple Nash equilibria of the followers' mean-field games, directly analyzing the game value of the Stackelberg mean-field game is formidable. Instead, we reformulate the Stackelberg mean-field game as a minimax optimization problem. This reformulation is inspired by the recent work on the equivalence between mean-field games and constrained optimization problems in \cite{guo2022mf}, and by incorporating appropiately the role of the leader. 
 
This new analytical framework for Stackelberg mean-field games enables us 
to characterize analytically the set of Nash equilibria and avoids the restrictive uniqueness  assumption of  Nash equilibrium in existing approaches. In particularly, when the leader's action space is finite, finding the leader's best action is reduced to solving finitely many optimization problems. 
When the leader's action space is not finite, the problem can be rewritten explicitly as a max-min optimization problem, for which computational methods such as gradient descent ascent \cite{lin2020gradient,jin2020local,lin2020near} and extragradient algorithms \cite{lee2021fast,chavdarova2019reducing,gidel2018variational} can be adopted. 

Furthermore, the characterization of Nash equilibria  facilitates analyzing the regularity  and the robustness properties of the game value. 
We show that the game value is not continuous, as feasible regions of constrained optimization problems can be sensitive with respect to slight modification of constraints. Moreover, when there is model uncertainty, the game value for the leader suffers non-vanishing sub-optimality. In order for her to obtain a near-optimal solution, she needs to be more pessimistic with anticipation of model errors and adopts a relaxed version of the original Stackelberg game. Finally, the optimization framework developed here can be generalized for other classes of Stackelberg games including the principal-agent problem.

\paragraph{Other related works.} 

Mean-field type Stackelberg game with  continuous-time models setting has also attracted lots of research interests recently. It stems from an important generalization of mean-field games to incorporate a major player who interacts with a large number of minor players. This problem is initially introduced in \cite{huang2010large} and \cite{nguyen2012linear} under linear-quadratic settings, and has been extended to more general settings \cite{nourian2013}, 
and cases with finite state \cite{carmona2016finite}. These initial works all focus on finding Nash equilibrium among all players without a hierarchical structure.  Unlike these, Stackelberg game under mean field setting have been studied extensively as well. \cite{bensoussan2016mean} starts with a general model setup and studies necessary conditions for the optimal controls and equilibrium condition, with their follow up work \cite{bensoussan2015mean} working with delayed information and \cite{bensoussan2017linear} finding explicit solution under the linear-quadratic framework. Along the same line, a local decentralized control obtained from a mean-field game is provided in \cite{moon2018linear} to construct an $\epsilon$-Stackelberg-Nash equilibrium for the corresponding $N$-player game under linear-quadratic framework; state feedback equilibrium is studied in \cite{huang2020mean} and the existence, the uniqueness and time-consistency of the Stackelberg mean-field  game as a bi-level Mckean-Vlasov forward-backward stochastic differential equations has been studied in \cite{fu2020mean}, \cite{huang2020mean}, \cite{yang2021linear}. 

Another line of the continuous-time setting is the contract theory as the principal-agent games where the principal proposes a contract
to incentivize the agent(s) to work for the principal and receive a reward.  The seminal paper \cite{sannikov2008continuous} develops the continuous-time method to study this type of problems. In \cite{elie2019tale}, the authors study the problem with one principal and an infinite number of agents, which has plenty of applications in different areas, including energy optimal demand-response (see \cite{elie2021mean}, \cite{campbell2021deep}, \cite{mastrolia2022agency}), and epidemic containment (see \cite{aurell2022optimal}). More recently, the deep learning methods to solve a mean-field principal-agent games is also presented in \cite{campbell2021deep}. 

However, all above works analyze Stackelberg games assuming that the followers/agents will achieve certain specific equilibrium. As mentioned earlier, the direction of our work is different: it is to study the Stackelberg games when the leader may not have control over the followers and would optimize the worst-case cost, with focus on regularity and sensitivity analysis.

\paragraph{Organization of the paper}
The rest of the paper is organized as follows. Section \ref{sec:framework} presents the mathematical framework of the Stackelberg mean-field game. Section \ref{sec:smfg-optimization} shows how Stackelberg mean-field games are transformed into optimization problems. Section \ref{sec:sa_ep} discusses the regularity and sensitivity properties of Stackelberg mean-field games. In particular, Section \ref{subsec:discont} shows the discontinuity and instability of the game values with respect to $\epsilon$ and transitions/rewards and Section \ref{subsec:cont} establishes the robustness of relaxed Stackelberg mean-field games. The proofs of the results are in Section \ref{Sec:proof}, followed by a discussion on the connections between Stackelberg games and principal-agent problems.

\section{Mathematical Framework of Stackelberg Mean-Field Game}\label{sec:framework}
In this section, we introduce the mathematical framework for a discrete-time Stackelberg mean-field game in a finite-time horizon $T$. 
\subsection{Game setup}
In this game, there is one leader and an infinite number of identical,  indistinguishable, and rational followers. The game starts by the leader taking an action  $a^\l$ from an action set $\A^\l$. This choice of $a^\l\in \A^\l$ by the leader is critical for the Stackelberg game, as will be clear shortly. Given the leader's  action $a^\l$,  the game of the followers is assumed to be a mean-field type, to be defined precisely below.

We emphasize that  the leader takes an action only at the beginning of the game as in an ad auction, where the platform (i.e., the leader) may fix one specific strategy or design for weeks or months, while the advertisers (i.e., the followers) may participate in the game whenever there is a slot to bid.

Throughout the paper, we will denote $\A^\l, \A^f$  respectively for the action sets of the leader and of the followers,  $\S^\l, \S^f$  for their respective state spaces, and $\PP(E)$ for the set of all distributions on a general space $E$. We will  assume that $|\S^f|,|\S^\l|$, $|\A^f|$ are finite. 


\paragraph{Follower's mean-field game given the leader's action $a^\l$.}
 
 
In the mean-field game of the followers,  at any  time $t=0,1,\dots,T-1$ when the  representative follower is at state $s_t^f\in \S^f$, he takes an action $a^f_t\sim \pi_t(\cdot|s_t^f)$, with $\pi_t:\S^f\to\PP(\A^f)$  a randomized policy  of the follower at time $t$. Given his state $s_t^f$, his action $a_t^f$, and the population state-action distribution of all followers denoted as {$L^f_t\in\PP(\S^f\times\A^f)$} at time $t$, he will receive a reward $r^{f,a^\l}_t(s_t^f,a_t^f,L_t^f)$, with $r_t^{f,a^\l}:\S^f\times\A^f\times\PP(\S^f\times\A^f)\to \R$, and move to the next state $s_{t+1}^f$ according to a transition probability $\P^{f,a^\l}_t(\cdot|s_t^f,a_t^f,L_t^f)$, where  $\P_t^{f,a^\l}:\S^f\times\A^f\times\PP(\S^f\times\A^f)\to \PP(\S^f)$.

In this game, the representative follower interacts with the population distribution flow $\bm L^f=(L^f_t)_{t=0}^T\in\PP(\S^f\times\A^f)^{T+1}$, and aims at finding the best policy to maximize his total expected rewards. For any policy $\bm\pi=\big\{\pi_t(\cdot|s^f):\S^f\to\PP(\A^f),\;t\in\{0,..,T\}\big\}$,  his  objective function {$J^{f,a^\l}$}, given  $\bm\pi$ and $\bm L^f=(L^f_t)_{t=0}^T$, is 
\begin{equation*}
J^{f,a^\l}(\bm\pi,\bm L^f):= \E^{f,a^\l,\bm\pi,\bm L^f}\bigg[\sum_{t=0}^T r_t^{f,a^\l}(s^f_t, a^f_t, L^f_t)\bigg], 
\end{equation*}
and his value function  $V^{f,a^\l}:\PP(\S^f\times\A^f)^{T+1}\to \R$ is
\begin{equation*}
    V^{f,a^\l}(\bm L^f)=\sup_{\bm\pi} J^{f,a^\l}(\bm\pi,\bm L^f),
\end{equation*}
where $\P^{f,a^\l,\bm\pi,\bm L^f}$ and $\E^{f,a^\l,\bm\pi,\bm L^f}$ denote respectively as the probability and expectation under this probability distribution for the representative follower when $s_0^f\in\mu^f_0,s_{t+1}^f\sim\P^{f,a^\l}_t(\cdot|s_t^f,a_t^f,L_t^f)$ for $t\in\{0,..,T-1\}$ and $a^f_t\sim \pi_t(\cdot|s_t^f)$ for $t\in\{0,..,T\}$.

We will denote for notation simplicity  $\M^{f,a^\l}:=(T,\S^f,\A^f,\mu_0^f,\P^{f,a^\l},r^{f,a^\l})$ for this mean-field game of followers, given the leader's action $a^\l$.

\paragraph{Leader's game.}\label{sec:leader}

Meanwhile, at time $t=0,\dots,T$, given the leader's action $a^\l$ (taken at the start of the game), her current state $s_t^\l$, and the population distribution flow of the followers $L^f_t$, she will receive a reward $r^\l_t(s_t^\l,a^\l,L^f_t):\S^\l\times\A^\l\times\PP(\S^f\times\A^f)\to \R$, and move to a new state $s^\l_{t+1}$ according to a probability distribution $\P^\l_t(\cdot|s_t^\l,a^\l,L^f_t)$, with $\P_t^{\l}:\S^\l\times\A^f\times\PP(\S^f\times\A^f)\to \PP(\S^\l)$. 
The leader aims at finding the best action that maximizes her total expected reward function $R^\l:\A^\l\times\PP(\S^f\times\A^f)^{T+1}\to \R$ given the population state-action distribution $\bm L$,  defined as 
\begin{equation*}
    R^\l(a^\l,\bm L^f):= \E^{\l,a^\l,\bm L^f}\bigg[\sum_{t=0}^T r_t^\l(s^\l_t, a^\l, L_t^f)\bigg].
\end{equation*}
Here $\P^{\l,a^\l,\bm L^f},\E^{\l,a^\l,\bm L^f}$ denote respectively as the probability and expectation  under this probability distribution for the leader when $s_0^\l\sim\mu^\l_0,s_{t+1}^\l\sim\P^\l_t(\cdot|s_t^\l,a^\l,L^f_t)$ for $t\in\{0,..,T-1\}$.

\subsection{Game value}
To analyze this class of Stackelberg mean-field games, we first adopt the notion of $\epsilon$-Nash equilibrium for the mean-field game of followers. 

\begin{Definition}[$\epsilon$-Nash equilibrium of mean-field games] \label{def:epnash}
A pair $(\bm\pi, \bm L^f)$ is an $\epsilon$-Nash equilibrium for the mean-field game $\M^{f,a^\l}=(T,\S^f,\A^f,\mu_0^f,\P^{f,a^\l},r^{f,a^\l})$ if
\begin{itemize}
    \item Policy $\bm\pi$ is $\epsilon$-suboptimal for any follower under the population distribution flow $\bm{L^f}$, \textit{i.e.},
    \[
    J^{f,a^\l}(\bm\pi,\bm L^f) \geq V^{f,a^\l}(\bm L^f)-\epsilon.
    \]
    \item Population distribution flow $\bm L^f$ is consistent  when all followers take policy $\bm\pi$, in the sense that
    \[
    \P^{f,a^\l,\bm\pi,\bm L^f}(s^f_t=s^f,a^f_t=a^f)=L^f_t(s^f,a^f),\qquad\text{for all }t\in\{0,..,T\},\;s^f\in\S^f,\;a^f\in\S^f.
    \]
\end{itemize}
\end{Definition}

This relaxed notion where followers aim at finding $\epsilon$-Nash equilibria instead of Nash equilibria captures sub-optimal strategies in games and is related to the bounded rationality in \cite{aumann1997rationality}, first proposed by Herbert A. Simon in \cite{simon1955behavioral} in economic theory and has since been widely adopted in psychology \cite{kahneman2003maps}, law \cite{korobkin2003bounded}, political science \cite{jones1999bounded}, and cognitive science \cite{lewis2014computational}. 
When $\epsilon=0$, it is the classical Nash equilibrium where the followers have perfect rationality; when $\epsilon>0$, the quantity of $\epsilon$ characterizes the tolerance of sub-optimality of the followers and therefore the distance to perfect rationality.

Next, we elaborate on the choice of $a^\l$ for the leader, which is crucial for defining the value of  this Stackelberg mean-field game. 
 
\paragraph{Leader's choice of action.} 
Recall that in the Stackelberg mean-field game, the leader and the followers may have conflicting interests and the followers may adopt different Nash equilibria from those desired by the leader. We therefore postulate that the objective of the leader is to find a best action $a^{\l,\star}$ that maximizes her total expected reward over the worst-case scenario over the $\epsilon$-best response set  $\text{BR}^\epsilon(a^\l)$, defined by 
\begin{equation*}
    \text{BR}^\epsilon (a^\l):=\Big\{(\bm\pi,\bm L^f): (\bm\pi,\bm L^f)\text{ is an }\epsilon\text{-Nash equilibrium for the game } \M^{f,a^\l}\Big\}.
\end{equation*}
The corresponding objective function $J^\l_\epsilon:\A^\l\to \R$ for the leader is then
\begin{equation}\label{eqn:valueleader}
    J^\l_\epsilon(a^\l)=\min_{(\bm\pi,\bm L^f)\in \text{BR}^\epsilon(a^\l)}R^\l(a^\l,\bm L^f).
\end{equation}

Finally,  with the $\epsilon$-Nash equilibria and  the $\epsilon$-best response set $\text{BR}^\epsilon(a^\l)$ of the followers,   the  corresponding  $(\epsilon,\epsilon')$-Stackelberg Nash equilibrium (SNE) for the Stackelberg mean-field game can be defined as follows.
\begin{Definition}[$(\epsilon,\epsilon')$-Stackelberg Nash equilibrium for the Stackelberg mean-field game]\label{defn:SNESMFG}
A pair $(\bm\pi^\star,\bm L^{f,\star},a^{\l,\star})$ is an $(\epsilon,\epsilon')$-Stackelberg Nash equilibrium if
\begin{itemize}
    \item $(\bm\pi^\star,\bm L^{f,\star})$ is an $\epsilon$-Nash equilibrium under leader's action $a^{\l}$, that is
    \begin{equation*}
        (\bm\pi^\star,\bm L^{f,\star})\in{\rm{BR}}^\epsilon(a^{\l,
        \star}).
    \end{equation*}
    \item Action $a^{\l,\star}$ is $\epsilon'$-optimal to the leader, such that
    \[
    J_\epsilon^\l(a^{\l,\star})>V^\l(\epsilon)-\epsilon'.
    \]
    where the value function of the leader is defined by
    \begin{equation*}
    V^\l(\epsilon):=\sup_{a^\l}J^{\l}_{\epsilon}(a^\l).
\end{equation*}
\end{itemize}
\end{Definition}

This notion  of $(\epsilon,\epsilon')$-Stackelberg Nash equilibrium is well defined as its existence is directly implied by the existence of Nash equilibrium of mean-field games (see Proposition 1 in \cite{guo2022mf}). For the purpose of completeness, we rephrase the results by the following proposition. 
\begin{Proposition}[Exsitence of Stackelberg Nash equilibrium]
    Suppose $\P_t^{f,a^\l}(\tilde s^f|s^f,a^f,L^f)$ and $r_t^{f,a^\l}(s^f,\\a^f,L^f)$ are both continuous in $L^f\in\PP(\S^f\times\A^f)$ for any $a^\l\in\mathcal{A}^\l$, $\tilde s^f,s^f\in\S^f$, $a^f\in\A^f$ and $t\in\{0,..,T\}$. Then for any $\epsilon\geq0$ and any $\epsilon'>0$, there exists an $(\epsilon,\epsilon')$-Stackelberg Nash equilibrium $(\bm\pi,\bm L^f,a^\l)$ for $\epsilon\geq0,\epsilon'>0$. In particular, 
    when $\A^\l$ is finite, the existence holds for $\epsilon'=0$.
\end{Proposition}

\subsection{Example} We end this section with an example of Stackelberg mean-field games with multiple Nash-equilibria and hence $\epsilon$-Nash equilibria among the followers, to further explain  the mathematical setup.  
 This example is adapted from \cite{guo2022mf}. 
 
\begin{Example}[Following the majority]\label{example:0}
Consider $\mathcal{S}^f=\mathcal{A}^f=\{1,\dots,n\}$ for the followers, with a deterministic transition probability of the following form: 
\begin{equation*}
    \P_t^f(j'|i,j,L_t^f)=\mathbbm{1}_{\{j=j'\}}, \quad \forall i,j,j'=1,\dots,n,\quad t=0,\dots,T.
\end{equation*}
The followers receive no reward at $t=0$. At $t>0$, if the population gather at the  state $i$, the representative follower will only receive a positive reward $r^i$ when he is also at state $i$ and chooses to stay at state $i$. If the population state-action joint distribution deviates from the state and action of the representative follower, his reward will be adjusted according to the distance between him and the followers, defined by
\begin{itemize}
    \item $r^f_0(i,j,L^f_0)=0$, \quad$\forall i,j=1,\dots,n$;
    \item $r^f_t(i,j,L^f_t)=\mathbbm{1}_{\{i=j\}}(r^i-r^i\|L^f_t-\bm e_{i,i}\|_2^2/2)$,\quad$\forall i,j=1,\dots,n,\quad t=1,\dots,T$;
\end{itemize}
where $r^i>0$ for all $i=1,\dots,n$ and $\bm e_{i,i}$ denote the unit vector in $\R^{|\S^f|\times |\A^f|}$ whose support is $(i,i)$. Such transition probabilities and reward functions incentivize all followers to gather at the same state and never leave. 

In this game, one can check that for any given positive vectors $\bm r=(r^1,\dots,r^n)$, there are at least $n$ different Nash equilibria. That is, for any $k=1,\dots,n$, let
\begin{itemize}
    \item $\pi_t(j|i)=\mathbbm{1}_{\{j=k\}}$  $\forall i,j=1,\dots,n$, \quad $t=0,\dots,T$;
    \item $L^f_0(i,j)=\mu_0(i)\pi_0(j|i)$ and $L_t^f(i,j)=\mathbbm{1}_{\{i=k,j=k\}}$, $\forall i,j=1,\dots,n$ and $t=1,\dots,T$,
\end{itemize}
then $(\bm\pi,\bm L^f)$ is a Nash equilibrium for the followers.

In the particular case of $n=2$, 
the leader's reward depends on $L^f_t$ and $a^\l=(r^1,r^2)$ in the following form:
\begin{equation}
    r^\l_t(s_t^\l,a^\l,L_t^f)=(r^2-r^1)\mathbbm{1}_{\{L_t^f=\mathbf{1}_{\{s=1,a=1\}}\}}+(r^1-r^2)\mathbf{1}_{\{L_t^f=\mathbbm{1}_{\{s=2,a=2\}}\}},
\end{equation}
which means that if all followers gather at state $1$, the leader will obtain a reward of $r^2-r^1$; if all followers gather at state $2$, the leader will get a reward of $r^1-r^2$; otherwise, the leader will get zero reward. 
Therefore, the action set of the leader in this game  effectively consists of choosing $r^1-r^2$; and the leader and the followers may have opposite goals. For instance, when $r^1>r^2$, the leader would prefer the followers to be at state $2$; yet the followers may favor state $1$ which gives a better reward $r^1$. 

If the leader  has the full control of the followers, she can first dictate the followers to their Nash equilibrium state $1$, and then set $r^1$ and $r^2$ so that $r^1-r^2$  is as large as possible. 
However, if the  leader does not have control over the followers, then she will consider the worst case scenario where the followers may gather at the Nash equilibrium state $2$. Thus the best action for the leader, according to the $(\epsilon, \epsilon')$-Stackelberg Nash equilibrium defined in \ref{defn:SNESMFG}, is to set $r^1=r^2$ to guarantee a non-negative and in fact a zero reward. Consequently,by considering the worst-case scenario among all possible Nash equilibria, she avoids getting a negative rewards in case  the followers gather at their instead of her favorite state. 
\end{Example}

\section{Analyzing Stackelberg mean-field games as Constrained Optimization Problems}\label{sec:smfg-optimization}
To analyze  this Stackelberg mean-field games in the sense of Definition \ref{defn:SNESMFG}, it is to solve for 
 the objective function of the leader $J^\l_\epsilon(a^\l)$: 
\begin{align}\label{RSMFG}
J^\l_\epsilon(a^\l)&=\text{minimize}_{(\bm\pi,\bm L^f)}  \quad \mathbb{E}^{\l,a^\l,\bm L^f}\bigg[\sum_{t=0}^T r^\l_t(s_t^\l, a^\l,L_t^f)\bigg] \\
&\text{subject to} \quad (\bm\pi,\bm L^f)\in\text{BR}^\epsilon(a^\l).\label{RSMFG:cons}
\end{align}
However, it is generally difficult to directly characterize the set of best responses ${\rm{BR}}^\epsilon(a^\l)$. Instead, we will consider the following constrained optimization problem:
\begin{align}
     \min_{(\bm \mu^\l,\bm d^f)} \quad  \sum_{t=0}^T\sum_{s^\l\in\S^\l}  \mu_t^\l(s^\l)r^\l_t(s^\l,a^\l,d^f_t)&  \label{op:begin}  \\
    \text{subject to} \quad  \sum_{s^\l\in\S^\l}\mu_t^\l(s^\l)\P_t^\l(\tilde s^\l|s^\l,a^\l,d^f_t)&=\mu_{t+1}^\l(\tilde s^{\l}),\quad\ \forall \tilde s^\l\in\S^\l,\;t\in\{0,..,T-1\},\label{eq:cons_occ_l}\\
    \sum_{s^f\in\S^f}\sum_{a^f\in A^f}d^f_t(s^f,a^f)\P^{f,a^\l}_t(\tilde s^f|s^f,a^f,d^f_t)& =\sum_{a^f\in\A^f}d^f_{t+1}(\tilde s^f,a^f), \quad\forall \tilde s^f\in\S^f,t\in\{0,..,T-1\},\label{eq:cons_occ_f}\\
     \sum_{a^f\in\A^f}d^f_t(s,a)&=\mu_0^f(s^f),\qquad\qquad\qquad\forall s^f\in\S^f,\label{eq:cons_init_f}\\
    \sum_{t\in \{0,..,T\}}\sum_{s^f\in\S^f}\sum_{a^f\in\A^f}  d^f_t(s^f,a^f)r^{f,a^\l}_t(s^f,a^f,d^f_t)\ &\geq V^{f,a^\l}(\bm d^f)-\epsilon,\label{eq:sub_opt_f}\\
     \mu_t^\l(s^\l)&\geq 0,\qquad\qquad\ \forall s^\l\in\S^\l,t\in\{0,..,T\},\\
     d^f_t(s^f,a^f)&\geq 0,\qquad\qquad\ \forall s^f\in\S^f,a^f\in\A^f,t\in\{0,..,T\}.\label{op:end}
\end{align}
Here  $V^{f,a^\l}(\bm d^f)$ is the  value function of the representative follower, given the leaders action $a^\l$ and mean-field flow $\bm d^f$. 

This optimization problem is derived by introducing occupation measure $\bm d^f$ for the representative follower in the mean-field game, as in \cite{guo2022mf}, and by characterizing the two conditions of $\epsilon$-Nash equilibrium using the occupation measure: first \eqref{eq:cons_occ_f} and \eqref{eq:cons_init_f} for the consistency condition and then \eqref{eq:sub_opt_f} for the $\epsilon$-suboptimality condition. Note that the occupation measure of the representative follower should be consistent with the mean-field flow.

We will establish the following equivalence.

\begin{Proposition}\label{prop:opwithvalue}
    The optimization problem \eqref{RSMFG}-\eqref{RSMFG:cons} is equivalent to the optimization problem \eqref{op:begin}-\eqref{op:end} in the following sense: 
    \begin{itemize}
        \item if $(\bm \pi,\bm d^f)$ achieves an optimal value of \eqref{RSMFG}-\eqref{RSMFG:cons}, then $(\bm\mu^\l,\bm d^f)$ achieves an optimal value of the optimization problem \eqref{op:begin}-\eqref{op:end} where
        \begin{equation}\label{def:dmul}
        \mu^\l_t(s^\l):=\P^{\l,a^\l,\bm d^f}(s^\l_t=s^\l).
        \end{equation}
        \item if $(\bm\mu^\l,\bm d^f)$ achieves an optimal value of the optimization problem \eqref{op:begin}-\eqref{op:end}, then $(\bm\pi,\bm d^f)$ achieves an optimal value of \eqref{RSMFG}-\eqref{RSMFG:cons}, where 
        \begin{equation}\label{def:pimuf}
            \begin{cases}
            \pi_t(a^f|s^f)=\frac{d^f_t(s^f,a^f)}{\sum_{\tilde a^f\in\A^f}d^f_t(s^f,\tilde a^f)},&\text{when }\sum_{\tilde a^f\in\A^f}d^f_t(s^f,\tilde a^f)>0.
            \\
            \pi_t(\cdot|s^f)\text{ being any probability vector},&\text{when }\sum_{\tilde a^f\in\A^f}d^f_t(s^f,\tilde a^f)=0.
            \end{cases}
        \end{equation}
    \end{itemize}
\end{Proposition}

Furthermore, due to the fact that when the leader's action $a^\l$ and the mean-field flow $\bm d^f$ are fixed, the optimization problem for the representative follower becomes an MDP problem, and its value $V^{f,a^\l}(\bm d^f)$ 
can  be represented explicitly via  the linear program formulation of MDPs and its associated KKT conditions.

\begin{Proposition}\label{prop:valuelp}
The  value function of the representative follower, $V^{f,a^\l}(\bm d^f)$ , given the leader's action $a^\l$ and occupation measure $\bm d^f$, has the representation of   
    \begin{align}
        &V^{f,a^\l}(\bm d^f)=-\bm c^{a^\l}(\bm d^f)^\top \bm x,\label{eq:kktbegin}\\
        \text{where }\quad&\bm A^{a^\l}(\bm d^f) \bm x=\bm b(\mu_0^f),\quad \bm A^{a^\l}(\bm d^f)^\top \bm u+\bm v=\bm c^{a^\l}(\bm d^f),\label{eq:kktmid}\\
        & \bm v^\top \bm x=0,\quad \bm x\geq 0,\quad \bm v\geq 0,\label{eq:kktend}
    \end{align}
for some $\bm x, \bm u, \bm v$. Here $\bm c^{a^\l}(\bm d^f)\in \mathbb{R}^{|\S^f||\A^f|(T+1)}$, $\bm b(\mu_0^f)\in\mathbb{R}^{|\S^f|(T+1)}$ are defined by 
    \begin{equation}\label{kkt_param1}
    \bm c^{a^\l}(\bm d^f):=\left[
\begin{array}{c}
-r^{f,a^\l}_0(\cdot,\cdot,d_0^f)\\
\vdots\\
-r^{f,a^\l}_T(\cdot,\cdot,d_T^f)
\end{array}
\right],\quad
\bm b(\mu_0^f):=\left[
\begin{array}{c}
0\\
\vdots\\
0\\
\mu_0^f
\end{array}
\right],
    \end{equation}
and $\bm A^{a^\l}(\bm d^f)\in\mathbb{R}^{|\S^f|(T+1)\times|\S^f||\A^f|(T+1)}$ is defined as: 
\begin{equation}\label{kkt_param2}
\bm A^{a^\l}(\bm d^f)=\left[
\begin{array}{ccccccc}
   W_0(d_0^f)  & -Z &  0 & 0 & \cdots &0 & 0\\
   0 & W_1(d_1^f)  & -Z &  0 & \cdots &0 & 0\\ 
   0 & 0 & W_2(d^f_2) & -Z & \cdots & 0 &0 \\
   \vdots & \vdots & \vdots & \vdots & \ddots & \vdots &\vdots \\
   0 & 0 & 0 & 0 & \cdots & W_{T-1}(d^f_{T-1}) & -Z\\
   Z & 0 & 0 & 0 & \cdots & 0 & 0
\end{array}\right],
\end{equation}
where $W_t(d^f_t)\in \mathbb{R}^{|\S^f|\times |\S^f||\A^f|}$ is the matrix with the $l$-th row ($l=1,\dots,S$) being the  flattened vector  $[\P^{f,a^\l}_t(l|\cdot,\cdot,d^f_t)]\in\mathbb{R}^{SA}$ (with column-major order),  
and the matrix $Z$ is defined as 
\begin{equation}\label{eq:defn-Z}
 Z:=[\overbrace{I_{S},\dots,I_{S}]}^{|\A^f|}\in \mathbb{R}^{|\S^f|\times |\S^f||\A^f|},  
\end{equation}
where $I_S$ is the identity matrix with dimension $|\S^f|$.

\end{Proposition}

Combining Propositions \ref{prop:opwithvalue} and \ref{prop:valuelp}, we have the following.
\begin{Theorem}\label{thm:opfram}
Solving   problem \eqref{RSMFG}-\eqref{RSMFG:cons} is equivalent to solving the following optimization problem
\begin{align}
    \min_{(\bm \mu^\l,\bm d^f,\bm x,\bm u,\bm v,V)} \quad  \sum_{t=0}^T\sum_{s^\l\in\S^\l} \mu_t^\l(s^\l)r^\l_t(s^\l,a^\l&,  d^f_t).    \label{fixal:begin} \\
    \mathrm{subject}\;\mathrm{to} \quad  \sum_{s^\l\in\S^\l}\mu_t^\l(s^\l)\P_t^\l(\tilde s^\l|s^\l,a^\l,d^f_t)&=\mu_{t+1}^\l(\tilde s^{\l}),\qquad\qquad\forall \tilde s^\l\in\S^\l,\;t\in\{0,..,T-1\},\label{fixal:consist1}\\
    \sum_{s^f\in\S^f}\sum_{a^f\in A^f}d^f_t(s^f,a^f)\P^{f,a^\l}_t(\tilde s^f|s^f,a^f,d^f_t)&=\sum_{a^f\in\A^f}d^f_{t+1}(\tilde s^f,a^f), \quad\;\forall \tilde s^f\in\S^f,t\in\{0,..,T-1\},\\
     \sum_{a^f\in\A^f}d^f_t(s,a)&=\mu_0^f(s^f)\;\qquad\qquad\qquad\forall s^f\in\S^f,\label{fixal:consist3}\\
    \sum_{t\in \{0,..,T\}}\sum_{s^f\in\S^f}\sum_{a^f\in\A^f} d^f_t(s^f,a^f)r^{f,a^\l}_t(s^f,a^f,d^f_t)&\geq V-\epsilon,\label{fixal:ep1}\\
    V=-\bm c^{a^\l}(\bm d^f)^\top \bm x,\quad,\bm A^{a^\l}(\bm d^f)\bm x&=\bm b(\mu_0^f),\quad \bm A^{a^\l}(\bm d^f)^\top \bm u+\bm v=\bm c^{a^\l}(\bm d^f),\label{fixal:ep2}\\
     \bm v^\top \bm x&=0,\quad \bm x\geq 0,\quad \bm v\geq 0,\label{fixal:ep3}
    \\
    \mu_t^\l(s^\l)&\geq 0,\qquad\qquad \forall s^\l\in\S^\l,t\in\{0,..,T\},\\
     d^f_t(s^f,a^f)&\geq 0,\qquad\qquad \forall s^f\in\S^f,a^f\in\A^f,t\in\{0,..,T\}.\label{fixal:end}
\end{align}

where $\bm A^{a^\l},\bm b,\bm c^{a^\l}$ are set as \eqref{kkt_param1} and \eqref{kkt_param2}. Specifically, 
\begin{itemize}
    \item if $(\bm \pi,\bm d^f)$ achieves an optimal value of \eqref{RSMFG}-\eqref{RSMFG:cons}, then there exist some $(\bm \mu^\l,\bm d^f,\bm x,\bm u,\bm v,V)$ which achieves the optimal value of the optimization problem \eqref{fixal:begin}-\eqref{fixal:end}. More precisely, one can define $(\bm\mu^\l,\bm d^f)$ as in \eqref{def:dmul}, $V=V^{f,a^\l}(\bm d^f)$, and $\bm x,\bm u,\bm v$ as the corresponding auxiliary variable and solve \eqref{eq:kktbegin}-\eqref{eq:kktend} to achieve the optimal value.
    \item if $(\bm\mu^\l,\bm d^f,\bm x,\bm u,\bm v,V)$ achieves the optimal value of the optimization problem \eqref{fixal:begin}-\eqref{fixal:end}, then $(\bm \pi,\bm d^f)$ defined by \eqref{def:pimuf} achieves optimal value of \eqref{RSMFG}-\eqref{RSMFG:cons}.
\end{itemize} 
\end{Theorem}

This equivalent formulation facilitates the  regularity and sensitivity analysis of the Stackelberg mean-field game, it enables characterizing $\epsilon$-Nash equilibria of the Stackelberg mean-field games through the feasible regions of the associated optimization problem.

\section{Regularity and Sensitivity  of Stackelberg Mean-Field Games}
\label{sec:sa_ep} 
In this section, we study the sensitivity and robustness properties of the Stackelberg mean-field game. 
 
Let us first define the perturbation of a Stackelberg mean-field game.

\begin{Definition}[$(\delta_p,\delta_r)$-perturbation of Stackelberg mean-field game $\M^\l$]\label{perturb_def}
Take a Stackelberg
 mean-field game  $\M^\l=(T,\S^\l,\S^f,\A^\l,\A^f,\mu_0^\l,\mu_0^f,\{\M^{f,a^\l}\}_{a^\l\in\A^\l},\P^\l,r^\l)$. Then a new game
$\hat\M^l=(T,\S^\l,\S^f,\A^\l,\\
\A^f,\hat \mu_0^\l,\hat \mu_0^f,\{\hat{\M}^{f,a^\l}\}_{a^\l\in\A^\l},\hat \P^\l,\hat r^\l)$ is called a $(\delta_p,\delta_r)$-perturbation of  $\M^\l$,  if for each $\tilde s^f,s^f,a^f,s^\l,a^\l,L^f,$
\begin{equation*}
    \begin{aligned}
    \big|\P_t^{f,a^\l}(\tilde s^f|s^f,a^f,L^f)-\hat\P_t^{f,a^\l}(\tilde s^f|s^f,a^f,L^f)\big|\leq \delta_p,
    \\
    \big|\P_t^{\l}(\tilde s^\l|s^\l,a^\l,L^f)-\hat \P_t^{\l}(\tilde s^\l|s^\l,a^\l, L^f)\big|\leq \delta_p,
    \\
    \big\|\mu^\l_0-\hat\mu^\l_0\big\|_\infty\leq \delta_p,\quad\big\|\mu^f_0-\hat \mu^f_0\big\|_\infty\leq \delta_p.
    \\
    \big|r_t^{f,a^\l}(s^f,a^f,L^f)-\hat r_t^{f,a^\l}(s^f,a^f,\tilde L^f)\big|\leq \delta_r,
    \\
    \big|r^\l_t(s^\l,a^\l,L^f)-\hat r_t^{\l}(s^\l,a^\l, L^f)\big|\leq \delta_r.
    \end{aligned}
\end{equation*}
\end{Definition}

Next, we will present various examples of Stackelberg mean-field games to show its discontinuity with respect to $\epsilon$ and its instability with respect to perturbations in reward functions and in transition probabilities.

Then, we will show that both the objective function and the value function of the leader are right continuous with respect to $\epsilon$, and that one can achieve the robustness of the game by appropriate choices of  $\epsilon$. 



\subsection{Discontinuity and Instability of the Game Value}\label{subsec:discont}
In this section, we show that the Stackelberg mean-field game is sensitive to the perturbations. Specifically, the value of the game $J^\l_\epsilon(a^\l)$ is not necessarily continuous with respect either to the value of $\epsilon$ or to the perturbations of the transitions and rewards in the game. In particular, there may be some non-vanishing gap when the leader applies an optimal action from the perturbed model to the true model. 

\subsubsection{Instability with respect to $\epsilon$}
The first example shows that the objective function of the leader is not necessarily continuous with respect to $\epsilon$. 
\begin{Example}\label{example:1}
Consider a predator-prey type of Stackelberg game: the prey have two locations to congregate, the first is state $e$, where they are exposed to the predator and where they need to stay close to fight against the predator; the reward function $\mu^f(e)-\epsilon_0$ represents the positive correlation with the population distribution at state $e$ denoted by $\mu^f(e)$. The second one is state $s$, which is a safe shelter and the prey always have a fixed reward $1$ at the safe shelter. Mathematically, the follower's game $\M^{f,a^\l}$ with time horizon $T=1$ has  the state and action spaces of $\S^f=\A^f=\{e,s\}$, where $e$ represents the vulnerable state and $s$ represents the safe state. The transition probability $\P^f_0(\tilde s^f|s^f,a^f,L^f) =\mathbbm{1}\{\tilde s^f=a^f\}$, the reward function $r_0^{f,a^\l}(s^f,a^f,L^f)=0$, and
\begin{equation*}
    r_1^{f,a^\l}(s^f,a^f,L^f) = \begin{cases}
    \mu^{f}(e)-\epsilon_0,&s^f=e,
    \\
    1,&s^f=s,
    \end{cases}
\end{equation*}
where $\mu^{f}(e):=L^f(e,s)+L^f(e,e)$ is the marginal distribution for the state under $L^f$. 

Clearly, all prey choosing to staying at location $s$ at step $1$ is a Nash equilibrium (and therefore an $\epsilon$-Nash equilibrium) for the prey.
We will see that for a given $\epsilon$, when $\epsilon_0$ is sufficiently small such that $\epsilon\ge \epsilon_0$, 
then state $e$ is also an $\epsilon$-Nash equilibrium. 

To see this, take  a policy $\bm\pi$, then
\begin{equation*}
    \P^{f,\bm \pi,\bm L^f}(s_1^f=s^f) = \E^{f,\bm\pi,\bm L^f}[\pi_0(s^f|s_0^f)],
\end{equation*}
and the objective function $J(\bm\pi,\bm L^f)$ is
\begin{equation*}
    \begin{aligned}
    J(\bm\pi,\bm L^f)&=\E^{f,\bm\pi,\bm L^f}\big[r_1^f(s_1^f,a_1^f,L_1^f)\big]= \E^{f,\bm\pi,\bm L^f}[\pi_1(s|s_0^f)]+\E^{f,\bm\pi,\bm L^f}[\pi_0(e|s_0^f)] (\mu_1^f(e)-\epsilon_0)
    \\
    &=1+ \E^{f,\bm \pi,\bm L^f}[\pi_1(e|s_0^f)] (\mu_1^f(e)-\epsilon_0-1).
    \end{aligned}
\end{equation*}
Therefore it is easy to see that $V^f(\bm L^f) = 1$ for all $\bm L^f$. 

Now recall the definition of $\epsilon$-Nash equilibrium,
\begin{equation*}
    \begin{cases}
        J(\bm \pi,\bm L^f)\geq 1-\epsilon,
        \\
        \mu_1^f(e)=\E^{f,\bm \pi,\bm d_1^f}[\pi_1(e|s_0^f)],
    \end{cases}
\end{equation*}
we know that $(\bm\pi,\bm L^f)$ is an $\epsilon$-Nash equilibrium if and only if 
\begin{equation*}
\mu^f_1(e)(\mu^f_1(e)-\epsilon_0-1)+ \epsilon\geq 0.
\end{equation*} 

Now solving this inequality $\mu^f_1(e)(\mu^f_1(e)-\epsilon_0-1)+ \epsilon\geq 0$, we get
\begin{equation*}
    \mu_1^f(e)\leq \frac{(1+\epsilon_0)-\sqrt{(1+\epsilon_0)^2-4\epsilon}}2,\;\;\text{or }\mu_1^f(e)\geq \frac{(1+\epsilon_0)+\sqrt{(1+\epsilon_0)^2-4\epsilon}}2.
\end{equation*}
Note that $\mu_1^f(e)\leq 1$, the second inequality is valid if and only if $\epsilon\geq\epsilon_0$.

Meanwhile, let us consider a specific predator's action $a^\l$: hunting at location $e$. This action provides the predator with a reward $1-\mu^f(e)$, which is negatively related to the gathering rate of the prey as the higher number of crowd the predator faces, the more difficult the hunting is:  
\begin{equation*}
    r_0^\l(s^\l,a^\l,\bm L^f)=0,\;r_1^\l(s^\l,a^\l,\bm L^f)=r_1^\l(\bm L^f):=1-\mu^f_1(e),
\end{equation*}
then by the definition of the leader's objective function,
\begin{equation*}
    J_\epsilon^\l(a^\l)=\min_{(\bm\pi,\bm d^f)\in\text{BR}^\epsilon}(1-\mu_1^f(e))=\begin{cases}
    1-\frac{(1+\epsilon_0)-\sqrt{(1+\epsilon_0)^2-4\epsilon}}2, & 0\leq\epsilon<\epsilon_0,
    \\
    0, & \epsilon\geq \epsilon_0,
    \end{cases}
\end{equation*}
which is not continuous with respect to $\epsilon$ at $\epsilon_0$.
\\
\vspace{10pt}
\end{Example}

\subsubsection{Instability with respect to transitions and rewards}
The sensitivity of the Stackelberg mean-field game can also be seen through perturbing the game in Example \ref{example:1}. More precisely, if one considers  $\hat\M^\l$, a perturbed version of the Stackelberg mean-field game $\M^\l$ with perturbed objective function $\hat J_\epsilon^\l(a^\l)$, then $\hat J_\epsilon^\l(a^\l)$ does not converge to $J_\epsilon^\l(a^\l)$ even as the perturbation decreases to $0$. In the following,  we will use $\hat J^{\l,\delta_p,\delta_r}_\epsilon(a^\l)$ to highlight the dependency of the objective function on the perturbations $\delta_p,\delta_r$.

\begin{Example} [Example \ref{example:1} continued]\label{example:2}
Take Example \ref{example:1} and specify its parameter $\epsilon_0$ to be the same as the tolerance parameter $\epsilon$, then
\begin{equation*}
    J^\l_\epsilon(a^\l)=0.
\end{equation*}
Now, assume there exists a judgement error of $\delta_r$, meaning that reward function for the prey at state $e$ is changed to $\mu^f(e)-\epsilon_0-\delta_r$. That is, define $\hat\P_0^{f,a^\l},\hat r^{f,a^\l}_0,\hat r^\l= \P_0^{f,a^\l},r_0^{f,a^\l},r^\l$ and let
\begin{equation*}
     \begin{aligned}
        \hat r_1^{f,a^\l}(s^f,a^f, L^f) = \begin{cases}
        \mu^f(e)-\epsilon_0-\delta_r,&s^f=e,
        \\
        1,&s^f=s.
        \end{cases}
    \end{aligned}
\end{equation*}
Then this new game is a $(\delta_p,\delta_r)$-perturbation of the original game in Example \ref{example:1} for all $\delta_p>0$ by definition \ref{perturb_def}. Now the condition for the $\epsilon$-Nash equilibrium becomes
\begin{equation*}
    \hat\mu_1^f(e)\big(\hat\mu_1^f(e)-\epsilon_0-\delta_r-1\big)+\epsilon\geq 0,
\end{equation*}
and the assumption of $\epsilon_0=\epsilon$ implies
\begin{equation*}
    \hat\mu_1^f(e)\leq \frac{(1+\epsilon_0+\delta_r)-\sqrt{(1+\epsilon_0+\delta_r)^2-4\epsilon_0}}2.
\end{equation*}
Therefore, 
\begin{equation*}
    \hat J^{\l,\delta_p,\delta_r}_\epsilon(a^\l)=1-\frac{(1+\epsilon_0+\delta_r)-\sqrt{(1+\epsilon_0+\delta_r)^2-4\epsilon_0}}2.
\end{equation*}
And the objective function of the perturbed model for the predator does not converge to the true objective function even if the perturbation $\delta_r \to 0$:
\begin{equation*}
    \inf_{\delta_p,\delta_r\in(0,1]} \hat J^{\l,\delta_p,\delta_r}_\epsilon(a^\l)=1-\sup_{\delta_p,\delta_r\in(0,1]}\frac{(1+\epsilon_0+\delta_r)-\sqrt{(1+\epsilon_0+\delta_r)^2-4\epsilon_0}}2=1-\epsilon_0>0=J^\l_\epsilon(a^\l).
\end{equation*}

\end{Example}

Moreover, this gap in value does not vanish even if the leader is to choose his optimal action  
    $a^{\l,\delta_p,\delta_r,*}:=\argmax_{a^\l} \hat J^{\l,\delta_p,\delta_r}_\epsilon (a^\l).$
Here $\{a^{\l,\delta_p,\delta_r,*}\}_{\delta_p, \delta_r}$ is sub-optimal in the sense that there exists an $\eta>0$ such that
\begin{equation*}
    \sup_{\delta_p,\delta_r\in(0,1]} J_\epsilon^{\l}(a^{\l,\delta_p,\delta_r,*})\leq V^\l(\epsilon)-\eta.
\end{equation*}

\begin{Example}[Example \ref{example:1} continued]\label{example:3}
Now consider another perturbed version of Example \ref{example:1}. Let $\A^\l=\{g,l\}$, where $g$ represents the predator going out for hunting with full efforts and $l$ represents the predator being lazier. When the predator chooses action $l$, it will get a discounted reward of $\frac{1}{3}(1-\mu_1^f(e))+\frac{1-\epsilon_0}{3}$. Keeping  $\hat\P_0^{f,a^\l},\hat r^{f,a^\l}_0,\hat r^{f,a^\l}_0,\hat r_0^\l$ the same as in Example \ref{example:2} and $\P_0^{f,a^\l},r^{f,a^\l}_0,r^{f,a^\l}_0,r_0^\l$ the same as in Example \ref{example:1} when $\epsilon_0=\epsilon$ and
\begin{equation*}
    \hat r_1^\l(s^\l,a^\l,L^f)= \begin{cases}
    1-\mu_1^f(e),&a^\l=g,
    \\
    \frac{1}{3}(1-\mu_1^f(e))+\frac{1-\epsilon_0}{3},&a^\l=l.
    \end{cases}
\end{equation*}
Clearly this new model is a $(\delta_p,\delta_r)$-perturbation of the original model in Example \ref{example:1} for all $\delta_r, \delta_p>0$ by definition \ref{perturb_def}.
Now the claim is that the predator will choose action $g$ in the perturbed game: according to the previous calculation, $\hat J^{\l,\delta_p,\delta_r}_\epsilon(g)\geq 1-\epsilon_0$ and 
\begin{equation*}
\hat J_\epsilon^{\l,\delta_p,\delta_r}(l)=\frac{1}{3}\hat J_\epsilon^{\l,\delta_p,\delta_r}(g)+\frac{1-\epsilon_0}{3}<\hat J_\epsilon^{\l,\delta_p,\delta_r}(g).
\end{equation*}
Therefore, $a^{\l,\delta_p,\delta_r,*}=g$. Meanwhile, it is straightforward to see that $J^\l_\epsilon(g)=0$, $J^\l_\epsilon(l)=\frac{1-\epsilon_0}3$, which leads to
\begin{equation*}
    V^\l(\epsilon) = \frac{1-\epsilon_0}3 > \sup_{\delta_p,\delta_r\in(0,1]}J_\epsilon^\l(a^{\l,\delta_p,\delta_r,*}),
\end{equation*}
again a non-trivial value gap from the optimal choice of the predator. 
\end{Example}

\paragraph{Discussions.}  
As we can see in the above Example \ref{example:3}, a savvy predator may avoid the nontrivial value gap 
by considering optimizing over all $(\epsilon+\epsilon')$-Nash equilibria for some $\epsilon'>0$. This relaxation incorporates the case of $\mu^f(e)=1$ for the prey and enables the predator to regain its optimal action of being ``lazier''. In a sense, the predator should be more pessimistic when it anticipates model error. This is a key intuitive insight that leads to the robustness and stability properties for the objective functions, discussed in the next section.

\subsection{Continuity of Stackelberg mean-field game and robustness of  relaxed Stackelberg mean-field game}\label{subsec:cont}
Despite the discontinuity and instability of the Stackelberg mean-field game, one can show that the value of Stackelberg mean-field game is right continuous with respect to $\epsilon.$

\begin{Proposition}\label{prop_cont0}
 The objective function of the leader given her action $a^\l$, $J_\epsilon^\l(a^\l)$ as defined in \eqref{eqn:valueleader} is right  continuous with respect to $\epsilon$.
\end{Proposition}

In fact, this right continuity property for the game value is preserved for the value function that we recall
\begin{equation*}
    V^\l(\epsilon)=\sup_{a^\l}J^{\l}_{\epsilon}(a^\l).
\end{equation*}
\begin{Proposition}\label{prop_contmax}
    $V^\l(\epsilon)$ is right continuous with respect to $\epsilon$.
\end{Proposition}

We will postpone the proofs of these two propositions in Section \ref{Sec:proof}. The right continuity property enables us to show that with appropriate Lischitz assumptions, one can bound the gap between the value of the original Stackelberg mean-field game $\M^\l$ and that of its perturbed version $\hat\M^{\l}$. In particular, if one considers an appropriate relaxed version of the original $\M^\l$, then the optimal solution in the sense of Stackelberg Nash equilibrium for  this relaxed version leads to an $(\epsilon,\epsilon')$-Stackelberg Nash equilibrium
for $\M^\l$. 

To this end, we will need the assumption that the transitions and rewards for the leader and the followers are all Lipschitz continuous with respect to the mean-field flow.
\begin{Assumption}[Lipschitz conditions]\label{assump:lip}
 There exists a constant $C$ independent of the choice of $t,\tilde s^f,s^f,a^f,s^\l,a^\l$, such that for any $L^f, \tilde{L}^f\in\PP(\S^f\times\A^f)$,
\begin{equation*}
    \begin{aligned}
    \big|\P_t^{f,a^\l}(\tilde s^f|s^f,a^f,L^f)-\P_t^{f,a^\l}(\tilde s^f|s^f,a^f,\tilde L^f)\big|\leq C\|L^f-\tilde L^f\|_\infty,
    \\
    \big|r_t^{f,a^\l}(s^f,a^f,L^f)-r_t^{f,a^\l}(s^f,a^f,\tilde L^f)\big|\leq C\|\L^f-\tilde L^f\|_\infty,
    \\
    \big|\P_t^{\l}(\tilde s^\l|s^\l,a^\l,L^f)-\P_t^{\l}(\tilde s^\l|s^\l,a^\l,\tilde L^f)\big|\leq C\|L^f-\tilde L^f\|_\infty,
    \\
    \big|r^\l_t(s^\l,a^\l,L^f)-r_t^{\l}(s^\l,a^\l,\tilde L^f)\big|\leq C\|L^f-\tilde L^f\|_\infty.
    \end{aligned}
\end{equation*}
\end{Assumption}

The following lemma bounds the difference between value functions of the representative follower and the distributions under different models with the same policy, which is essential for the main results. 

\begin{Lemma}\label{lem_closeJ}
Assume that model $\M^\l$  satisfies Assumption \ref{assump:lip}. Let $\hat \M^\l$ be any $(\delta_p,\delta_r)$-perturbation of $\M^\l$. Denote $V^{f,a^\l}(\bm d^f)$ and $\hat V^{f,a^\l}(\hat{\bm d}^f)$ the value functions of the representative follower under the population distribution flows $\bm d^f,\hat {\bm d}^f$ under models $\M^\l$, $\hat \M^\l$ and the corresponding occupation measure for the leader being $\mu^\l(s^\l)$ and $\hat\mu^\l(s^\l)$ respectively. Assume further that $\bm d^f,\hat {\bm d}^f$ follow the same policy in that  there exists $\{(\pi_t)(a^f|s^f)\}_{t,a^f,s^f}$ such that
\begin{equation*}
   d^f_t(s^f,a^f)=\mu^f_t(s^f)\pi_t(a^f|s^f),\quad \hat d^f_t(s^f,a^f)=\hat \mu^f_t(s^f)\pi_t(a^f|s^f)
\end{equation*}
 where $\mu^f_t(s^f)=\sum_{a^f\in\A^f}d_t^f(s^f,a^f)$ and $\hat\mu^f_t(s^f)=\sum_{a^f\in\A^f}\hat d_t^f(s^f,a^f)$ for all $s^f\in\S^f, t\in\{0,..,T\}$. Then
\begin{equation*}
    \begin{aligned}
    \max\Big\{\big|\mu^\l_t(s^\l)-\hat \mu^\l_t(s^\l)\big|,\big|{\mu_t^f}(s^f)-\hat \mu_t^f(s^f)\big|,\big|d_t^f(s^f,a^f)-\hat d_t^f(s^f,a^f)\big|\Big\}\leq (C+S+1)^{t}\delta_p.
    \end{aligned}
\end{equation*}
Moreover, 
\begin{equation*}
    \big|V^{f,a^\l}(\bm d^f)-\hat V^{f,a^\l}(\hat{\bm d}^f)\big|\leq (C+S)(1+C+S)^{T+2}\delta_p+(T+1)\delta_r.
\end{equation*}
\end{Lemma}

\begin{Theorem}\label{thm:close}
Consider an $\M^\l$  satisfying Assumption \ref{assump:lip}, take $\hat \M^\l$  any $(\delta_p,\delta_r)$-perturbation of $\M^\l$, set $S=\max\{|\S^\l|,|\S^f|\}$ and  any $\epsilon'>0$ such that 
\begin{equation}\label{inequa:deltaepsilon}
    (C+S)(1+C+S)^{T+1}\delta_p+(T+1)\delta_r\leq\frac {\epsilon'}{2}.
\end{equation}
 Then for any $\epsilon\geq 0$,
\begin{equation}\label{inequ:big}
    \begin{aligned}
    J^\l_\epsilon(a^\l)&\geq \hat J^\l_{\epsilon+\epsilon'}(a^\l)-(C+S)(1+C+S)^{T+1}\delta_p-(T+1)\delta_r,
    \end{aligned}
\end{equation}
and if $\epsilon\geq\epsilon'$, 
\begin{equation}\label{inequ:small}
    J^\l_\epsilon(a^\l)\leq \hat J^\l_{\epsilon-\epsilon'}(a^\l)+(C+S)(1+C+S)^{T+1}\delta_p+(T+1)\delta_r.
\end{equation}
Moreover, for any $\delta>0$, there exists $\epsilon^*>0$, such that for any $0<\epsilon'<\epsilon^*$ and $\hat \M^\l$ $(\delta_p,\delta_r)$-close estimator of $\M^\l$ with $\delta_p,\delta_r$ satisfying \eqref{inequa:deltaepsilon} and $\hat a^{\l,\delta_p,\delta_r,\epsilon'}:=\argmax_{a^\l} \hat V^\l_{\epsilon+\epsilon'} (a^\l)$, we have 
\begin{equation*}
    J^\l_\epsilon (\hat a^{\l,\delta_p,\delta_r,\epsilon'})\geq V^\l(\epsilon) -\delta.
\end{equation*}
That is $a^{\l,\delta_p,\delta_r,\epsilon'}$ and its optimal value $(\bm\pi,\bm L^f)$ of \eqref{RSMFG}-\eqref{RSMFG:cons} is an $(\epsilon,\delta )$-Stackelbetrg Nash equilibrium under $\M^\l$.
\end{Theorem}

An immediate result is the following limiting property, which states that if the perturbation of the model can be arbitrarily small, one can obtain an approximate solution sufficiently close to the optimal one, as long as one enlarges the value of $\epsilon$ appropriately. We have shown in Section \ref{subsec:discont} that this limiting property does not hold for a fixed $\epsilon$.
\begin{Corollary}\label{coro:limit}
Consider an $\M^\l$  satisfying Assumption \ref{assump:lip}, set $S=\max\{|\S^\l|,|\S^f|\}$. Define the set 
\begin{equation}\label{ass-set}
E=\Big\{(\delta_p,\delta_r,\epsilon')\geq 0: (C+S)(1+C+S)^{T+1}\delta_p+(T+1)\delta_r\leq\frac {\epsilon'}{2}\Big\}.
\end{equation}
Then 
\begin{equation}
\lim_{\substack{\epsilon'\rightarrow 0\\ (\delta_p,\delta_r,\epsilon')\in E}} J^\l_\epsilon (\hat a^{\l,\delta_p,\delta_r,\epsilon'})=V^\l(\epsilon).
\end{equation}
where $\hat a^{\l,\delta_p,\delta_r,\epsilon'}$ is the optimal solution to the relaxed perturbed model for any $\hat \M^\l$  $(\delta_p,\delta_r)$-perturbation of $\M^\l$, that is 
\begin{equation*}
     \hat a^{\l,\delta_p,\delta_r,\epsilon'}:=\argmax_{a^\l} \hat J^\l_{\epsilon+\epsilon'} (a^\l),
\end{equation*}
and $\hat J^\l_{\epsilon+\epsilon'}(a^\l)$ is the value function under the relaxed perturbed model.

\end{Corollary}
\begin{Remark}
In both Theorem \ref{thm:close} and Corollary \ref{coro:limit}, we assume that the supreme of $\hat J^\l_{\epsilon+\epsilon'} (a^\l)$ can be achieved by some $\hat a^{\l,\delta_p,\delta_r,\epsilon}$ which is always true when the leader has only finitely many choices of action. When the leader has infinite choices, one may still reach the same result by carefully choosing the `optimal' action that is sufficiently close to the supreme.
\end{Remark}

\section{Proofs for the Main Results}\label{Sec:proof}
This section is devoted to the proofs of the main results. 

\paragraph{Proof of Proposition \ref{prop:opwithvalue}.}
The proof consists of two steps, the first  step is for the equivalence in terms of the feasibility and the second step for the equivalence in terms of the optimality.

\noindent\textbf{Step 1.} Suppose that $(\bm \pi,\bm d^f)$ satisfies \eqref{RSMFG:cons}, we will show that $(\bm \mu^\l,\bm d^f)$ defined by \eqref{def:dmul} satisfies  the constraints \eqref{eq:cons_occ_l}-\eqref{op:end}.

First, \eqref{eq:cons_occ_l} follows from the definition of $\P^{\l,a^\l,\bm d^f}$: \begin{align*}
    \mu^\l_{t+1}(\tilde s^\l)&=\P^{\l,a^\l,\bm d^f}(s^\l_{t+1}=\tilde s^\l)=\sum_{s^\l\in\S^\l}\P^{\l,a^\l,\bm d^f}(s^\l_{t+1}=\tilde s^\l,s^\l_t=s^\l)\\
    &=\sum_{s^\l\in\S^\l}\P^{\l,a^\l,\bm d^f}(s^\l_t=s^\l)\P^\l_t(\tilde s^\l|s^\l,a^\l,d_t^f)=\sum_{s^l\in\S^\l}\mu^\l_t(s^\l)\P^\l_t(\tilde s^\l|s^\l,a^\l,d^f_t).
\end{align*}
Secondly, \eqref{eq:cons_occ_f} holds because of the consistency condition in Definition \ref{def:epnash}:
\begin{equation*}
    \begin{aligned}
    \sum_{a^f\in\A^f}d_{t+1}^f(\tilde s^f,a^f)=&\mathbb{P}^{f,a^\l,\bm \pi,\bm d^f}(s^f_{t+1}=\tilde s^f)=\sum_{s^f\in\S^f}\sum_{a^f\in\A^f}\mathbb{P}^{f,a^\l,\bm \pi,\bm d^f}(s^f_{t+1}=\tilde s^f,s^f_{t}=s^f,a_t^f=a^f)
    \\
    =&\sum_{s^f\in\S^f}\sum_{a^f\in\A^f}\mathbb{P}^{f,a^\l,\bm \pi,\bm d^f}(s^f_{t}=s^f,a^f=a^f)\P^{f,a^\l}_t(\tilde s^f|s^f,a^f,\mu^f_t)
    \\
    =&\sum_{s^f\in\S^f}\sum_{a^f\in\A^f} d_t(s^f,a^f)\P^{f,a^\l}_t(\tilde s^f|s^f,a^f,\mu^f_t).
    \end{aligned}
\end{equation*}
Next,  the initial constraint\eqref{eq:cons_init_f} is due to 
\begin{equation*}
    \mu^f_0(s^f)=\mathbb{P}^{f,a^\l,\bm \pi,\bm d^f}(s_0^f=s^f)=\sum_{a^f\in\A^f}\mathbb{P}^{f,a^\l,\bm \pi,\bm d^f}(s_0^f=s^f,a_0^f=a^f)=\sum_{a^f\in\A^f}d_0^f(s^f,a^f).
\end{equation*}

As for \eqref{eq:sub_opt_f},  we have
\begin{equation*}
    J^{f,a^\l}(\bm \pi,\bm d^f)=\sum_{t=0}^T\sum_{s^f\in\S^f}\sum_{a^f\in\A^f} d^f_t(s^f,a^f)r^{f,a^\l}_t(s^f,a^f,d^f_t)
\end{equation*}
from 
$d^f_t(s^f,a^f)=\mathbb{P}^{f,a^\l,\bm \pi,\bm d^f}(s^f_t=s^f,a^f_t=a^f)$.
By the $\epsilon$-optimal condition in Definition \ref{def:epnash}, $J^{f,a^\l}(\bm \pi,\bm d^f)\\\geq V^{f,a^\l}(\bm d^f)-\epsilon$, hence \eqref{eq:sub_opt_f}. 
Finally, the last two conditions follow directly by definitions of $\bm \mu^\l$ and $\bm d^f$.

Conversely,  suppose that $(\bm \mu^\l,\bm d^f)$ satisfies  the constraints \eqref{eq:cons_occ_l}-\eqref{op:end}, we will show that $(\bm \pi,\bm d^f)$ defined by \eqref{def:pimuf} satisfies \eqref{RSMFG:cons}, i.e., $(\bm \pi,\bm d^f)$ is an $\epsilon$-Nash equilibrium under the model $\M^{f,a^\l}$. 
It suffices  to check the two conditions in Definition \ref{def:epnash}. For the consistency condition of the population distribution, for $t=1,2,..,T,$  define
\begin{equation*}
    \mu_t^f(s^f):=\sum_{a^f\in\A^f}d_t^f(s^f,a^f).
\end{equation*}
Note that by \eqref{eq:cons_init_f}, the above equation also holds for $t=0$. Now by induction, \eqref{eq:cons_occ_f} and $d^f_t(s^f,a^f)=\pi_t(a^f|s^f)\mu_t^f(s^f)$, we get
\begin{equation*}
    \mu_t^f(s^f)=\mathbb{P}^{f,a^\l,\bm \pi,\bm d^f}(s_t^f=s^f),\;\;d_t^f(s^f,a^f)=\mathbb{P}^{f,a^\l,\bm \pi,\bm d^f}(s_t^f=s^f,a_t^f=a^f),
\end{equation*}
hence the consistency of the population distribution.

To show the $\epsilon$-optimal condition, since $d^f_t(s^f,a^f)=\mu^f_t(s^f)\pi_t(a^f|s^f)=\P^{f,\pi,\bm d^f,a^\l}(s^f_t=s^f,a^f_t=a^f)$,  we have by \eqref{eq:sub_opt_f},
\begin{equation*}
   J^{f,a^\l}(\bm\pi,\bm d^f)=\sum_{t=0}^T\sum_{s^f\in\S^f}\sum_{a^f\in\A^f} d^f_t(s^f,a^f)r^{f,a^\l}_t(s^f,a^f,\mu^f_t)\geq V^{f,a^\l}(\mu^f)-\epsilon,
\end{equation*}
thus completing Step 1.

\noindent\textbf{Step 2.}
Now we  focus on the equivalence in terms of the optimality.

Suppose that $(\bm \pi,\bm d^f)$ achieves the optimal value of \eqref{RSMFG}-\eqref{RSMFG:cons}, we will show that $(\bm \mu^\l,\bm d^f)$ defined by \eqref{def:dmul} achieves the optimal value of the optimization problem \eqref{op:begin}-\eqref{op:end}.
Consider any feasible $(\tilde{\bm\mu}^{\l},\tilde {\bm d}^f)$ satisfying \eqref{eq:cons_occ_l}-\eqref{op:end}, one can construct its corresponding $(\tilde {\bm\pi},\tilde {\bm d}^f)$ by \eqref{def:pimuf}, which by Step 1 satisfies \eqref{RSMFG:cons}. Since $(\bm \pi,\bm d^f)$ achieves the optimal value of \eqref{RSMFG}-\eqref{RSMFG:cons},  $J^\l(a^\l,
\bm d^f)\leq J^\l(a^\l,\tilde{\bm d}^f)$. Then by the definition of $\bm \mu^\l$ and $\tilde{\bm\mu}^{\l}$,
\begin{equation*}
    \sum_{t=0}^T\sum_{s^\l\in\S^\l} \mu_t^\l(s^\l)r^\l_t(s^\l,a^\l,d^f_t)=J^\l(a^\l,\bm\mu^\l)\leq J^\l(a^\l,\tilde {\bm\mu}^\l)=\sum_{t=0}^T\sum_{s^\l\in\S^\l} \tilde\mu_t^\l(s^\l)r^\l_t(s^\l,a^\l,\tilde d^f_t).
\end{equation*}
 That is,  $(\bm \mu^\l,\bm d^f)$ achieves the optimal value of the optimization problem \eqref{op:begin}-\eqref{op:end}.

Conversely,  suppose that $(\bm \mu^\l,\bm d^f)$ achieves the optimal value of the optimization problem \eqref{op:begin}-\eqref{op:end}, we will show that $(\bm \pi,\bm d^f)$ defined by \eqref{def:pimuf} achieves the optimal value of \eqref{RSMFG}-\eqref{RSMFG:cons}. To see this, 
consider any $\epsilon$-Nash equilibrium $(\tilde{\bm\pi},\tilde{\bm d}^f)$, one can construct the corresponding tuple $(\tilde{\bm\mu}^{\l},\tilde{\bm d}^f)$ by \eqref{def:dmul}, which by Step 1 satisfies \eqref{eq:cons_occ_l}-\eqref{op:end}. By the optimality of $(\bm\mu^\l,\bm d^f)$,
\begin{equation*}
    \sum_{t=0}^T\sum_{s^\l\in\S^\l} \mu_t^\l(s)r^\l_t(s^\l,a^\l,d^f_t)\leq\sum_{t=0}^T\sum_{s^\l\in\S^\l}\tilde\mu_t^\l(s^\l)r^\l_t(s^\l,a^\l,\tilde d^f_t).
\end{equation*}
By \eqref{eq:cons_occ_l}, $\mu^\l_t(s^\l)=\P^{\l,a^\l,\bm d^f}(s^\l_t=s^\l)$ and $\tilde\mu^\l_t(s^\l)=\P^{\l,a^\l,\tilde{\bm d}^f}(s^\l_t=s^\l)$. Thus 
\begin{equation*}
    J^\l(a^\l,\bm d^f)=\sum_{t=0}^T\sum_{s^\l\in\S^\l} \mu_t^\l(s^\l)r^\l_t(s^\l,a^\l,d^f_t)\leq \sum_{t=0}^T\sum_{s^\l\in\S^\l} \tilde\mu_t^\l(s^\l)r^\l_t(s^\l,a^\l,\tilde d^f_t)=J^\l(a^\l,\tilde {\bm d}^f),
\end{equation*}
implying that $(\bm \pi,\bm d^f)$ achieves the optimal value of \eqref{RSMFG}-\eqref{RSMFG:cons}. \qed

\paragraph{Proof of Proposition \ref{prop:valuelp}.}
Following \cite[Lemma 2]{guo2022mf}, the value $V^{f,a^\l}(\bm d^f)$ equals the value of the following linear program:
\begin{equation}\label{primal}
    \begin{array}{ll}
        \text{minimize}_{\bm x} &-\bm c^{a^\l}(\bm d^f)^\top \bm x\\
        \text{subject to} &\bm A^{a^\l}(\bm d^f) \bm x=\bm b(\mu_0^f),\quad \bm x\geq0.
    \end{array}
    \end{equation}
By introducing dual variables $\bm u$ and $\bm v$, its corresponding dual problem is
\begin{equation}\label{dual}
    \begin{array}{ll}
        \text{maximize}_{\bm u, \bm v} & \bm b(\mu_0^f)^\top \bm u\\
        \text{subject to} &\bm A^{a^\l}(\bm d^f)^\top \bm u + \bm v=\bm c^{a^\l}(\bm d^f),\quad \bm v\geq0.
    \end{array}
    \end{equation}
The representation of $V^{f,a^\l}(\bm d^f)$ can be directly implied by the strong duality theorem and the KKT conditions \cite{luenberger1984linear}.  
\qed

\paragraph{Proof of Lemma \ref{lem_closeJ}.}
    Since $\hat{d}_t(s^f,a^f):=\hat\mu_t^f(s^f)\pi_t(a^f|s^f)$, $d_t(s^f,a^f):=\mu_t^f(s^f)\pi_t(a^f|s^f)$ and $\|\pi_t\|_\infty\leq 1$, it suffices to prove the result for $\bm\mu^f$ and $\bm\mu^\l$, by induction. When $t=0$, it holds by the definition of the $(\delta_p,\delta_r)$-close estimator. Suppose the inequality holds for $t$, by the definition of $\mu_{t+1}^\l$ and $\hat \mu_{t+1}^\l$, for any $\tilde s^\l\in\S^\l$, we have
    \begin{equation*}
        \begin{aligned}
        \big|\mu_{t+1}^\l(\tilde s^\l)-\hat \mu_{t+1}^\l(\tilde s^\l)\big| = \bigg|\sum_{s^\l\in\S^\l} \Big(\hat{\mu}_t^\l(s^\l)\hat\P_t^{\l}(\tilde s^\l|s^\l,a^\l,\hat{d}_t^f) - \mu_t^\l(s^\l)\P_t^{\l}(\tilde s^\l|s^\l,a^\l,{d}_t^f)\Big)\bigg|,
        \end{aligned}
    \end{equation*}
    and  each term inside the above summation can be rewritten as
    \begin{equation*}
        \begin{aligned}
        \Big(\hat{\mu}_t^\l(s^\l)-{\mu}_t^\l(s^\l)\Big)\hat\P_t^{\l}(\tilde s^\l|s^\l,\hat{d}_t^f) +{\mu}_t^\l(s^\l)\Big(&\big(\hat\P_t^{\l}(\tilde s^\l|s^\l,a^\l,\hat{d}_t^f)-\P_t^{\l}(\tilde s^\l|s^\l,a^\l,\hat{d}_t^f)\big)
        \\
        &+\big(\P_t^{\l}(\tilde s^\l|s^\l,a^\l,\hat{d}_t^f)-\P_t^{\l}(\tilde s^\l|s^\l,a^\l,{d}_t^f)\big)\Big).
        \end{aligned}
    \end{equation*}
    Splitting these three terms, by induction the first term can be bounded by
    \begin{equation*}
        \sum_{s^\l\in\S^\l}\big|\hat{\mu}_t^\l(s^\l)-{\mu}_t^\l(s^\l)\big|\hat\P_t^{\l}(\tilde s^\l|s^\l,a^\l,\hat{d}_t^f)\leq |\S^\l| (C+S+1)^{t}\delta_p\leq S(C+S+1)^{t}\delta_p.
    \end{equation*}
    By the definition of $(\delta_p,\delta_r)$-close estimator, the second term can be bounded by
    \begin{equation*}
         \sum_{s^\l\in\S^\l}{\mu}_t^\l(s^\l)\Big|\hat\P_t^{\l}(\tilde s^\l|s^\l,a^\l,\hat{d}_t^f)-\P_t^{\l}(\tilde s^\l|s^\l,a^\l,\hat{d}_t^f)\Big|\leq \sum_{s^\l\in\S^\l}{\mu}_t^\l(s^\l)\delta_p = \delta_p\leq (C+S+1)^{t}\delta_p.
    \end{equation*}
    By Assumption \ref{assump:lip}, the third term can be bounded by
    \begin{equation*}
         \sum_{s^\l\in\S^\l}{\mu}_t^\l(s^\l)\Big|\P_t^{\l}(\tilde s^\l|s^\l,a^\l,\hat{d}_t^f)-\P_t^{\l}(\tilde s^\l|s^\l,a^\l,{d}_t^f)\Big|\leq C\sum_{s^\l\in\S^\l}{\mu}_t^\l(s^\l)\|\hat d_t^f-d_t^f\|_\infty \leq C(C+S+1)^{t}\delta_p.
    \end{equation*}
    Combined, we get 
    \begin{equation*}
        \big|\mu_{t+1}^\l(\tilde s^\l)-\hat \mu_{t+1}^\l(\tilde s^\l)\big|\leq (C+S+1)^{t+1}\delta_p,
    \end{equation*}
    which completes the proof for $\mu^\l_t$. As for $\bm\mu^f_t$, the idea is exactly the same: write
    \begin{equation*}
         \big|\hat \mu_{t+1}^f(\tilde s^f)-\mu_{t+1}^f(\tilde s^f)\big|=\bigg|\sum_{s^f\in\S^f}\sum_{a^f\in\A^f} \big(\hat{d}_t(s^f,a^f)\hat\P_t^{f,a^\l}(\tilde s^f|s^f,a^f,\hat{d}_t^f)-{d}_t(s^f,a^f)\P_t^{f,a^\l}(\tilde s^f|s^f,a^f,{d}_t^f)\big)\bigg|,
    \end{equation*}
    and bound the first term by
    \begin{equation*}
        \begin{aligned}
        \sum_{s^f\in\S^f}&\sum_{a^f\in\A^f}\Big|\hat{d}_t^f(s^f,a^f)-{d}_t^f(s^f,a^f)\Big|\hat\P_t^{f,a^\l}(\tilde s^f|s^f,a^f,\hat{d}_t^f)
        \\
        &= \sum_{s^f\in\S^f}\big|\hat{\mu}_t^f(s^f)-{\mu}_t^f(s^f)\big|\sum_{a^f\in\A^f}\pi_t(a^f|s^f)\hat\P_t^{f,a^\l}(\tilde s^f|s^f,a^f,\hat{d}_t^f)
        \\
        &\leq |\S^f| (C+S+1)^{t}\delta_p\leq S(C+S+1)^{t}\delta_p,
        \end{aligned}
    \end{equation*}
    because $\sum_{a^f\in\A^f}\pi_t(a^f|s^f)\hat\P_t^{f,a^\l}(\tilde s^f|s^f,a^f,\hat{d}_t^f)=\hat\P^{f,a^\l,\bm\pi,\hat{\bm d}^f}(s_{t+1}^f=\tilde s^f|s^f_t=s^f)\leq 1$, and bound the second  and the third terms by
    \begin{equation*}
        \begin{aligned}
         \sum_{s^f\in\S^f}&\sum_{a^f\in\A^f}{d}_t^f(s^f,a^f)\Big|\hat\P_t^{f,a^\l}(\tilde s^f|s^f,a^f,\hat{d}_t^f)-\P_t^{f,a^\l}(\tilde s^f|s^f,a^f,\hat{d}_t^f)\Big|\leq (C+S+1)^{t}\delta_p.
         \\
         \sum_{s^f\in\S^f}&\sum_{a^f\in\A^f}{d}_t^f(s^f,a^f)\Big|\P_t^{f,a^\l}(\tilde s^f|s^f,a^f,\hat{d}_t^f)-\P_t^{f,a^\l}(\tilde s^f|s^f,a^f,{d}_t^f)\Big|\leq C(C+S+1)^{t}\delta_p,
        \end{aligned}
    \end{equation*}
    which completes the first part of the proof.
    
    The second part of the proof utilizes  the linear program formulation, i.e., \eqref{eq:kktbegin}-\eqref{eq:kktend}. Suppose that $\bm x$ achieves the optimal value for \eqref{eq:kktbegin}-\eqref{eq:kktend} under model $\M$, we can similarly define the policy as in \eqref{def:pimuf} such that
    \begin{equation*}
        \begin{cases}
        \tilde\pi_t(a^f|s^f)=\frac{x_t(s^f,a^f)}{\sum_{a^f\in\A^f}x_t(s^f,a^f)},&\text{when }\sum_{a^f\in\A^f}x_t(s^f,a^f)>0.
        \\
        \tilde\pi_t(\cdot|s^f)\text{ be any probability vector},&\text{when }\sum_{a^f\in\A^f}x_t(s^f,a^f)=0.
        \end{cases}
    \end{equation*}
    Now  define $\hat {\bm x}$ by 
    \begin{equation*}
        \begin{aligned}
        \hat{x}_0(s^f,a^f)&:=\hat \mu_0^f(s^f)\tilde\pi_0(a^f|s^f);
        \\
        \hat{x}_{t+1}(\tilde s^f,\tilde a^f)&:=\sum_{s^f\in\S^f}\sum_{a^f\in\A^f} \hat{x}_t(s^f,a^f)\hat\P_t^{f,a^\l}(\tilde s^f|s^f,a^f,\hat{d}_{t-1}^f)\tilde\pi_t(\tilde a^f|\tilde s^f),\;&\forall& t\in\{0,..,T-1\}.
        \end{aligned}
    \end{equation*}
    Then $\hat x$ satisfies the constraints \eqref{eq:kktbegin}-\eqref{eq:kktend} under model $\hat\M^\l$. Further, define
    \begin{equation*}
        y_t(s^f):=\sum_{a^f\in\A^f}x_t(s^f,a^f),\;\;\hat y_t(s^f)=\sum_{a^f\in\A^f}\hat x_t(s^f,a^f),
    \end{equation*}
    then by the same argument before, we have for all $s^f\in\S^f,\;a^f\in\A^f,\;t\in\{0,..,T\}$,
    \begin{equation*}
        \max\big\{|y_t(s^f)-\hat y_t(s^f)|,|x_t(s^f,a^f)-\hat x_t(s^f,a^f)|\big\}\leq (1+C+|\S^f|)^t\delta_p.
    \end{equation*}
    Next, let us consider the difference:
    \begin{equation*}
         \bigg|\sum_{t=0}^T\sum_{s^f\in\S^f}\sum_{a^f\in\A^f} \big(\hat x_t(s^f,a^f)\hat r^{f,a^\l}_t(s^f,a^f,\hat d^f_t)- x_t(s^f,a^f) r^{f,a^\l}_t(s^f,a^f,d^f_t)\big)\bigg|,
    \end{equation*}
    for which each term inside the summation can be written as the summation of three terms, with the first term  
    \begin{equation*}
        \begin{aligned}
        &\sum_{t=0}^T\sum_{s^f\in\S^f}\sum_{a^f\in\A^f} \big|x_t(s^f,a^f)-\hat x_t(s^f,a^f)\big|\cdot \big|r^{f,a^\l}_t(s^f,a^f,d^f_t)\big|
        \\
        \leq &\sum_{t=0}^T\sum_{s^f\in\S^f}\big|y_t(s^f)-\hat y_t(s^f)\big|\sum_{a^f\in\A^f} \pi_t(a^f|s^f) \big|r^{f,a^\l}_t(s^f,a^f,d^f_t)\big|
        \\
        \leq &|\S^f|\sum_{t=0}^T(1+C+|\S^f|)^t\delta_p\leq|\S^f| (1+C+|\S^f|)^T\delta_p
        \end{aligned}
    \end{equation*}
    since $\sum_{a^f\in\A^f} \pi_t(a^f|s^f) \big|r^{f}_t(s^f,a^f,\mu^f_t)\big|\leq 1$, the second term  
    \begin{equation*}
        \begin{aligned}
        &\sum_{t=0}^T\sum_{s^f\in\S^f}\sum_{a^f\in\A^f} x_t(s^f,a^f)\cdot \big|\hat r^{f,a^\l}_t(s^f,a^f,\hat d^f_t)-r^{f,a^\l}_t(s^f,a^f,\hat d^f_t)\big|
        \\
        \leq &\sum_{t=0}^T\sum_{s^f\in\S^f}\sum_{a^f\in\A^f}x_t(s^f,a^f)\delta_r=(T+1)\delta_r,
        \end{aligned}
    \end{equation*}
    and the third term  
    \begin{equation*}
        \begin{aligned}
        &\sum_{t=0}^T\sum_{s^f\in\S^f}\sum_{a^f\in\A^f} x_t(s^f,a^f)\cdot \big|r^{f,a^\l}_t(s^f,a^f,\hat d^f_t)-r^{f,a^\l}_t(s^f,a^f,d^f_t)\big|
        \\
        \leq &C\sum_{t=0}^T\sum_{s^f\in\S^f}\sum_{a^f\in\A^f}x_t(s^f,a^f)\|\hat d^f_t-d^f_t\|_\infty=C\sum_{t=0}^T\|\hat d^f_t-d^f_t\|_\infty\leq C(1+C+S)^{T+1}\delta.
        \end{aligned}
    \end{equation*}
    That is,
    \begin{equation*}
    \begin{aligned}
    \bigg|\sum_{t=0}^T\sum_{s^f\in\S^f}\sum_{a^f\in\A^f} \big(\hat x_t(s^f,a^f)\hat r^{f,a^\l}_t(s^f,a^f,\hat d^f_t)&- x_t(s^f,a^f) r^{f,a^\l}_t(s^f,a^f,d^f_t)\big)\bigg|
    \\
    &\leq (C+S)(1+C+S)^{T+1}\delta_p+(T+1)\delta_r,
    \end{aligned}
    \end{equation*}
    which implies that 
    \begin{equation*}
        \begin{aligned}
        \hat V^{f,a^\l}(\hat{\bm d}^f)&\geq \sum_{t=0}^T\sum_{s^f\in\S^f}\sum_{a^f\in\A^f} \big(\hat x_t(s^f,a^f)\hat r^{f,a^\l}_t(s^f,a^f,\hat d^f_t)
        \\
        &\geq \sum_{t=0}^T\sum_{s^f\in\S^f}\sum_{a^f\in\A^f} \big( x_t(s^f,a^f) r^{f,a^\l}_t(s^f,a^f, d^f_t)-(C+S)(1+C+S)^{T+1}\delta_p-(T+1)\delta_r
        \\
        &= V^{f,a^\l}(\bm d^f)-(C+S)(1+C+S)^{T+1}\delta_p-(T+1)\delta_r.
        \end{aligned}
    \end{equation*}
    Similarly,
    \begin{equation*}
        V^{f,a^\l}(\bm d^f)\geq \hat V^{f,a^\l}(\hat {\bm d}^f)-(C+S)(1+C+S)^{T+1}\delta_p-(T+1)\delta_r.
    \end{equation*} \qed

\paragraph{Proof of Theorem \ref{thm:close}.}
Let us first focus on \eqref{inequ:big}. Suppose $(\bm \mu^\l,\bm d^f,\bm x,\bm u,\bm v,V)$ achieves the  value of the optimization problem \eqref{fixal:begin}-\eqref{fixal:end}, one can construct $\bm\pi$ by \eqref{def:pimuf}. Similarly for \eqref{def:dmul}, define some $(\hat{\bm\mu}^\l,\hat{\bm d}^f)$ under the new model $\hat\M^\l$. That is, for all $\tilde s^f,s^f\in\S^f,a^f\in\A^f$ and $\tilde s^\l,s^\l\in\S^\l$ 
\begin{equation}\label{equ:hatmud}
    \begin{aligned}
    \hat{d}_t(s^f,a^f)&:=\hat\mu_t^f(s^f)\pi_t(a^f|s^f),\;&\forall& t\in\{0,..,T\},
    \\
    \hat \mu_{t+1}^f(\tilde s^f)&:=\P^{f,\bm\pi,\hat{\bm d}^f}(s_{t+1}^f=\tilde s^f)=\sum_{s^f\in\S^f}\sum_{a^f\in\A^f} \hat{d}_t(s^f,a^f)\hat\P_t^{f}(\tilde s^f|s^f,a^f,\hat{d}_t^f),\;&\forall& t\in\{0,..,T-1\},
    \\
    \hat \mu_{t+1}^\l(\tilde s^\l)&:=\P^{\l,\hat{\bm d}^f}(s_{t+1}^f=\tilde s^f)=\sum_{s^\l\in\S^\l} \hat{\mu}_t^\l(s^\l)\hat\P_t^{\l}(\tilde s^\l|s^\l,a^\l,\hat{d}_t^f),\;&\forall& t\in\{0,..,T-1\},
    \end{aligned}
\end{equation}
where $\hat{\bm\mu}^f:=\{\hat\mu_t^f\}_{t=1,2,..,T}$ being the auxiliary variables representing the marginal distribution for the state which also matches the definition of $\hat\mu_0^f$. By this construction, we know that $(\hat{\bm\mu}^\l,\hat{\bm d}^f)$ satisfies \eqref{fixal:consist1}-\eqref{fixal:consist3} for the model estimator $\hat M^\l$.

In order to prove \eqref{inequ:big}, we first show the existence of some $(\hat{\bm x},\hat{\bm u},\hat{\bm v},V)$ such that \eqref{fixal:ep1}-\eqref{fixal:ep3} holds for the model estimator $\hat M^\l$ with $\epsilon+\epsilon'$.

Let us define $\bm\mu_t^f:=\{\mu_t^f\}_{t=1,2,..,T}$, $\bm\hat\mu_t^f:=\{\hat \mu_t^f\}_{t=1,2,..,T}$ as the auxiliary variables representing the marginal distribution for the state by 
\begin{equation*}
    \mu_t^f(s^f)=\sum_{a^f\in\A^f}d^f(s^f,a^f), \quad \hat\mu_t^f(s^f)=\sum_{a^f\in\A^f}\hat d^f(s^f,a^f),
\end{equation*}
which matches the definition of $\mu_0^f$ and $\hat \mu_0^f$. Lemma \ref{lem_closeJ} provides the closeness of $(\bm\mu^\l,\bm\mu^f,\bm d^f)$ and $(\hat{\bm\mu}^\l,\hat{\bm\mu}^f,\hat{\bm d}^f)$

Now recall  Proposition \ref{prop:opwithvalue} and Proposition \ref{prop:valuelp} and consider $\hat V$ represents the value function for the followers given the population distribution $\hat{\bm \mu}$ and the corresponding $(\hat{\bm x},\hat{\bm y},\hat{\bm z})$ satisfy \eqref{fixal:ep2}-\eqref{fixal:ep3} for the model estimator $\hat M^\l$.

In order to check \eqref{fixal:ep1}, it suffices to show that 
\begin{equation*}
    \sum_{t\in \{0,..,T\}}\sum_{s^f\in\S^f}\sum_{a^f\in\A^f} \hat d_t(s^f,a^f)\hat r^{f,a^\l}_t(s^f,a^f,\hat d^f_t)\geq \hat V^{f,a^\l}(\hat {\bm d}^f)-\epsilon-\epsilon'. 
\end{equation*}
By the similar argument for Lemma \ref{lem_closeJ}, we have
\begin{equation*}
    \begin{aligned}
    \bigg|\sum_{t\in \{0,..,T\}}\sum_{s^f\in\S^f}&\sum_{a^f\in\A^f} \hat d_t(s^f,a^f)\hat r^{f,a^\l}_t(s^f,a^f,\hat d^f_t)-d_t(s^f,a^f)r^{f,a^\l}_t(s^f,a^f,d^f_t)\bigg|
    \\
    &\leq (C+S)(1+C+S)^{T+1}\delta_p+(T+1)\delta_r,
    \end{aligned}
\end{equation*}
Combining  Lemma \ref{lem_closeJ} and  \eqref{fixal:ep1}, we have 
\begin{equation*}
    \begin{aligned}
    \sum_{t\in \{0,..,T\}}&\sum_{s^f\in\S^f}\sum_{a^f\in\A^f} \hat d_t(s^f,a^f)\hat r^{f,a^\l}_t(s^f,a^f,\hat {d}^f_t)
    \\
    &\geq \sum_{t\in \{0,..,T\}}\sum_{s^f\in\S^f}\sum_{a^f\in\A^f} d_t(s^f,a^f)r^{f,a^\l}_t(s^f,a^f,d^f_t) - (C+S)(1+C+S)^{T+1}\delta_p-(T+1)\delta_r
    \\
    &\geq V^f(\bm d^f) - \epsilon - (C+S)(1+C+S)^{T+1}\delta_p-(T+1)\delta_r
    \\
    &\geq V^f(\hat {\bm d}^f)-\epsilon-2\big((C+S)(1+C+S)^{T+1}\delta_p+(T+1)\delta_r)
    \\
    &\geq V^f(\hat{\bm d}^f)-\epsilon-\epsilon',
    \end{aligned}
\end{equation*}
thus complete the proof of the feasibility of $(\hat {\bm\mu}^\l,\hat {\bm d}^f)$ of the optimization problem \eqref{fixal:begin}-\eqref{fixal:end}. 

Our next step is to bound
\begin{equation*}
    \bigg|\sum_{t=0}^T\sum_{s^\l\in\S^\l} \big(\hat \mu_t^\l(s^\l)\hat r^\l_t(s^\l,a^\l,\hat d^f_t)-\mu_t^\l(s^\l)r^\l_t(s^\l,a^\l,d^f_t)\big)\bigg|.
\end{equation*}
Similarly as in the proof of Lemma \ref{lem_closeJ}, we write out the three terms:
\begin{equation*}
    \begin{aligned}
    \sum_{t=0}^T\sum_{s^\l\in\S^\l} \big|\hat \mu_t^\l(s^\l)-\mu_t^\l(s^\l)\big|\cdot \big|\hat r^\l_t(s^\l,a^\l,\hat d^f_t)\big|&\leq& |\S^\l| (1+C+S)^{T+1}\delta_p,
    \\
    \sum_{t=0}^T\sum_{s^\l\in\S^\l} \mu_t^\l(s^\l) \big|\hat r^\l_t(s^\l,a^\l,\hat d^f_t)-r^\l_t(s^\l,a^\l,\hat d^f_t)\big|&\leq& (T+1)\delta_r,
    \\
    \sum_{t=0}^T\sum_{s^\l\in\S^\l} \mu_t^\l(s^\l) \big|r^\l_t(s^\l,a^\l,\hat d^f_t)-r^\l_t(s^\l,a^\l,d^f_t)\big|&\leq&C (1+C+S)^{T+1}\delta_p,
    \end{aligned}
\end{equation*}
and conclude that
\begin{equation*}
     \bigg|\sum_{t=0}^T\sum_{s^\l\in\S^\l} \big(\hat \mu_t^\l(s^\l)\hat r^\l_t(s^\l,a^\l,\hat  d^f_t)-\mu_t^\l(s^\l)r^\l_t(s^\l,a^\l,d^f_t)\big)\bigg|\leq(C+S)(1+C+S)^{T+1}\delta_p+(T+1)\delta_r.
\end{equation*}
Therefore, by Theorem \ref{thm:opfram},
\begin{equation*}
    \begin{aligned}
    \hat J^\l_{\epsilon+\epsilon'} (a^{\l})&\leq \sum_{t=0}^T\sum_{s^\l\in\S^\l} \hat \mu_t^\l(s^\l)\hat r^\l_t(s^\l,a^\l,\hat  d^f_t)
    \\
    &\leq\sum_{t=0}^T\sum_{s^\l\in\S^\l} \mu_t^\l(s^\l)r^\l_t(s^\l,a^\l,d^f_t)+(C+S)(1+C+S)^{T+1}\delta_p+(T+1)\delta_r
    \\
    &=J^\l_{\epsilon} (a^{\l})+(C+S)(1+C+S)^{T+1}\delta_p+(T+1)\delta_r.
    \end{aligned}
\end{equation*}
The other inequality can be deduced using exactly the same method by the optimality of $\hat J^\l_{\epsilon-\epsilon'}(a^\l)$ and by repeating the previous argument to construct a feasible solution under model $\M^\l$. Thus finish the proof of \eqref{inequ:big} and \eqref{inequ:small}.

To finish the proof, recall that by Proposition \ref{prop_cont0}, we can find $\epsilon^*>0$ such that $\epsilon^*\leq \frac\delta 4$ and for all $\epsilon'<\epsilon^*$
\begin{equation*}
    \sup_{a^\l}J^\l_{\epsilon+\epsilon'} (a^{\l})\geq V^\l(\epsilon) - \frac \delta 2.
\end{equation*}
Note that \eqref{inequ:big} and \eqref{inequ:small} lead to
\begin{equation*}
    \sup_{a^\l}J^\l_{\epsilon+\epsilon'} (a^{\l})\leq \sup_{a^\l}\hat J^\l_{\epsilon+\frac{\epsilon'}2} (a^{\l})+(C+S)(1+C+S)^{T+1}\delta_p+(T+1)\delta_r,
\end{equation*}
and 
\begin{equation*}
   J^\l_{\epsilon}(a^{\l,*})\geq \hat J^\l_{\epsilon+\frac{\epsilon'} 2}(a^{\l,*})-(C+S)(1+C+S)^{T+1}\delta_p-(T+1)\delta_r,
\end{equation*}
for all $\delta_p,\delta_r$ satisfying \eqref{inequa:deltaepsilon} and $a^{\l,*}:=\argmax_{a^\l} \hat J^\l_{\epsilon+\epsilon'} (a^\l)$ for any perturbed  $(\delta_p,\delta_r)$-estimator $\M^\l$. Putting these together, we get
\begin{equation*}
    \begin{aligned}
    J^\l_{\epsilon}(a^{\l,*})&\geq \sup_{a^\l}\hat J^\l_{\epsilon+\frac{\epsilon'}{2}} (a^{\l}) - (C+S)(1+C+S)^{T+1}\delta_p-2(T+1)\delta_r
    \\
    &\geq\sup_{a^\l}J^\l_{\epsilon+\epsilon'} (a^{\l}) - 2(C+S)(1+C+S)^{T+1}\delta_p-2(T+1)\delta_r
    \\
    &\geq\sup_{a^\l}J^\l_{\epsilon} (a^{\l}) - \frac\delta 2 - 2(C+S)(1+C+S)^{T+1}\delta_p-2(T+1)\delta_r
    \\
    &\geq\sup_{a^\l}J^\l_{\epsilon} (a^{\l}) - \frac\delta 2 - 2\epsilon' \geq V^\l(\epsilon) - \delta.
    \end{aligned}
\end{equation*} \qed

\paragraph{Proof of Proposition \ref{prop_cont0}.}
By the definition of the optimization problem,   $J_\epsilon^\l(a^\l)$ is decreasing with respect to $\epsilon$. It then suffices to show that $\liminf_{\epsilon\to \epsilon^*} J_\epsilon^\l(a^\l)\geq J_{\epsilon^*}^\l(a^\l)$ for all $\epsilon^*\geq 0$. 
We will prove this by contradiction. Suppose that $\liminf_{\epsilon\to \epsilon^*} J_\epsilon^\l(a^\l)<J_{\epsilon^*}^\l(a^\l)$, then by Proposition \ref{prop:opwithvalue}, there exists $\epsilon_n\to\epsilon^*$ and $(\bm\mu^{\l,n},\bm d^{f,n})$ such that $\liminf_{\epsilon\to \epsilon^*}J_\epsilon^\l(a^\l)=\lim_{n\to\infty} J_{\epsilon^n}^\l(a^\l)$ and $ J_{\epsilon^n}^\l(a^\l)=\sum_{t=0}^T\sum_{s^\l\in\S^\l} \mu_t^{\l,n}(s^\l)r^\l_t(s^\l,a^\l,d^{f,n}_t),$ where
\begin{equation*}
    \begin{aligned}
    & \mu_{t+1}^{\l,n}(\tilde s^{\l})=\sum_{s^\l\in\S^\l}\mu_t^{\l,n}(s^\l)\P_t^\l(\tilde s^\l|s^\l,a^\l,d^{f,n}_t),\qquad\qquad\qquad\;\;\;\forall \tilde s^\l\in\S^\l,\;t\in\{0,..,T-1\},\\
    &\sum_{a^f\in\A^f}d^{f,n}_{t+1}(s^f,a^f)=\sum_{s^f\in\S^f}\sum_{a^f\in \A^f}d^{f,n}_t(s^f,a^f)\P^{f,a^\l}_t(\tilde s^f|s^f,a^f,d^{f,n}_t),\quad \forall \tilde s^f\in\S^f,t\in\{0,..,T-1\},\\
    & \sum_{t\in \{0,..,T\}}\sum_{s^f\in\S^f}\sum_{a^f\in\A^f} d^n_t(s^f,a^f)r^{f,a^\l}_t(s^f,a^f,d^{f,n}_t)\geq V^{f,a^\l}(\bm d^{f,n})-\epsilon^n,
    \end{aligned}
\end{equation*}
where $\bm\mu^{\l,n}\in\PP(\S^\l)^{T},\bm d^{f,n}\in\PP(\S\times\A)^{T+1}$ are all within some compact set. Therefore, there exists a converging sub-sequence  denoted for simplicity still as $(\bm\mu^{\l,n},\bm d^{f,n})$  such that  
\begin{equation*}
    \big(\bm\mu^{\l,n},\bm d^{f,n}\big)\to \big(\bm\mu^{\l,*},\bm d^{f,*}\big).
\end{equation*}
By Assumption \ref{assump:lip},
\begin{equation*}
    \begin{aligned}
    \begin{aligned}
    & \mu_{t+1}^{\l,*}(\tilde s^{\l})=\sum_{s^\l\in\S^\l}\mu_t^{\l,*}(s^\l)\P_t^\l(\tilde s^\l|s^\l,a^\l,d^{f,*}_t),\qquad\qquad\qquad\;\;\;\forall \tilde s^\l\in\S^\l,\;t\in\{0,..,T-1\},\\
    &\sum_{a^f\in\A^f}d^{f,*}_{t+1}(s^f,a^f)=\sum_{s^f\in\S^f}\sum_{a^f\in \A^f}d^{f,*}_t(s^f,a^f)\P^{f,a^\l}_t(\tilde s^f|s^f,a^f,d^{f,*}_t),\quad \forall \tilde s^f\in\S^f,t\in\{0,..,T-1\},\\
    & \sum_{t\in \{0,..,T\}}\sum_{s^f\in\S^f}\sum_{a^f\in\A^f} d^{f,*}_t(s^f,a^f)r^{f,a^\l}_t(s^f,a^f,d^{f,*}_t)\geq \limsup_{n\to\infty} V^{f,a^\l}(\bm d^{f,n})-\epsilon^*,
    \end{aligned}.
    \end{aligned}
\end{equation*}
Moreover, using the same technique in the second part of the proof of Lemma \ref{lem_closeJ}, we can see $\limsup_{n\to\infty} V^{f,a^\l}(\bm d^{f,n})=V^{f,a^\l}(\bm d^{f,*})$, therefore $\big(\bm \mu^{\l,*},\bm d^{f,*}\big)$ is inside the feasible set which leads to 
\begin{equation*}
    \liminf_{\epsilon\to \epsilon^*} J_{\epsilon^n}^\l(a^\l)=\lim_{n\to\infty} J_{\epsilon^n}^\l(a^\l)=\sum_{t=0}^T\sum_{s^\l\in\S^\l} \mu_t^{\l,*}(s^\l)r^\l_t(s^\l,a^\l,d^{f,*}_t)\geq J_{\epsilon^*}^\l(a^\l),
\end{equation*}
a contradiction and proof completed.\qed

\paragraph{Proof of Proposition \ref{prop_contmax}.}
Fixed $\epsilon^*$ and for any $\delta>0$,  there exists some $a^{\l,*}$ such that 
\begin{equation*}
    V^\l(\epsilon^*)\leq J^\l_{\epsilon^*}(a^{\l,*})+\frac{\delta}{2},
\end{equation*}
then by the right continuity of $J^\l_\cdot(a^{\l,*})$ from Proposition \ref{prop_cont0},  there exists some $\epsilon_0>\epsilon^*$ such that for all $\epsilon^*\leq\epsilon\leq\epsilon_0$, we have
\begin{equation*}
    J^\l_{\epsilon^*}(a^{\l,*})\leq J^\l_{\epsilon}(a^{\l,*})+\frac{\delta}{2}. 
\end{equation*}
Thus, for all $\epsilon^*\leq\epsilon\leq\epsilon_0$,
\begin{equation*}
    V^\l(\epsilon^*) \leq J^\l_{\epsilon}(a^{\l,*})+\delta\leq  V^\l(\epsilon)+\delta,
\end{equation*}
which completes the proof, when combined with the fact that $V^\l(\cdot)$ is an decreasing function and $\delta$ is arbitrary.\qed

\section{Discussions}
Stackelberg games with potentially deviating followers 
are closely related to the  principal-agent problems, and  both are bi-level optimization problems.

Indeed, the principal's objective function $J^\l_\epsilon:\A^\l\to \R$ can be defined as
\begin{equation*}
    J^\l_\epsilon(a^\l)=\max_{(\bm\pi,\bm L^f)\in \text{BR}^\epsilon(a^\l)}J^\l(a^\l,\bm L^f),
\end{equation*}
since the principal can assign the agents the Nash equilibrium to maximize her objective and we can consider the action space of the leader as all measurable functions of the followers' states up to time $T$. Such definition of the objective function can also be applied to the Stackelberg game with an optimistic leader who aims at optimizing over the  best-case-scenario among all followers' possible responses.

Meanwhile, one can reformulate the principal agent problem of finding the best $a^\l$ as the following  optimization problem:
\begin{align*}
\text{maximize}_{(a^\l,\bm\pi,\bm L^f)}  \quad &\mathbb{E}^{\l,a^\l,\bm L^f}\bigg[\sum_{t=0}^T r^\l_t(s_t^\l, a^\l,L_t^f)\bigg], \\
\text{subject to} \quad &(\bm\pi,\bm L^f)\in\text{BR}^\epsilon(a^\l),
\end{align*}
using the same technique in this paper to deal with the set of best responses.
 
Note that instead of a minimax structure in the worst-case scenario for Stackelberg games considered in this paper, one has obtained a maximization optimization problem for the principal-agent which jointly optimizes $\bm\pi, \bm L^f$ and $a^\l$. 
Nevertheless, the instability property as in Examples \ref{example:1}, \ref{example:2} and \ref{example:3} for the Stackelberg game still holds for the principal-agent game, as can be seen by simply changing the reward functions of the leader to their negative values.

The sensitivity analysis of the Stackelberg games established in Section \ref{subsec:cont} does not directly apply to the principal-agent problems. When the action space for the leader is finite, Propositions \ref{prop_cont0} and \ref{prop_contmax} both hold and lead to the same conclusion as in Theorem \ref{thm:close}. However, when the action space for the leader is infinite, Proposition \ref{prop_cont0} is still true. But the value function for the leader $V^\l(\epsilon)$ may lose the right continuity property 
since $V^\l(\epsilon)$ becomes an increasing function and the supreme of a series of increasing right-continuous functions is not necessarily right-continuous. Thus, Proposition \ref{prop_contmax} may fail. In fact, sensitivity analysis for the principal-agent problem has not been studied in the existing literature and is worth future exploration.



\bibliographystyle{plain}
\bibliography{SMFG.bib}

\begin{thebibliography}{10}

\bibitem{alcantara2020repeated}
Guillermo Alcantara-Jim{\'e}nez and Julio~B Clempner.
\newblock Repeated {Stackelberg} security games: Learning with incomplete state
  information.
\newblock {\em Reliability Engineering \& System Safety}, 195:106695, 2020.

\bibitem{aumann1997rationality}
Robert~J Aumann.
\newblock Rationality and bounded rationality.
\newblock In {\em Cooperation: Game-Theoretic Approaches}, pages 219--231.
  Springer, 1997.

\bibitem{aurell2022optimal}
Alexander Aurell, Rene Carmona, Gokce Dayanikli, and Mathieu Lauriere.
\newblock Optimal incentives to mitigate epidemics: a {Stackelberg} mean field
  game approach.
\newblock {\em SIAM Journal on Control and Optimization}, 60(2):S294--S322,
  2022.

\bibitem{bai2021sample}
Yu~Bai, Chi Jin, Huan Wang, and Caiming Xiong.
\newblock Sample-efficient learning of {Stackelberg} equilibria in general-sum
  games.
\newblock {\em Advances in Neural Information Processing Systems}, 34, 2021.

\bibitem{bensoussan2017linear}
Alain Bensoussan, MHM Chau, Y~Lai, and Sheung Chi~Phillip Yam.
\newblock Linear-quadratic mean field stackelberg games with state and control
  delays.
\newblock {\em SIAM Journal on Control and Optimization}, 55(4):2748--2781,
  2017.

\bibitem{bensoussan2015mean}
Alain Bensoussan, Michael~HM Chau, and Sheung Chi~Phillip Yam.
\newblock Mean field {Stackelberg} games: Aggregation of delayed instructions.
\newblock {\em SIAM Journal on Control and Optimization}, 53(4):2237--2266,
  2015.

\bibitem{bensoussan2016mean}
Alain Bensoussan, Michael~HM Chau, and Sheung~CP Yam.
\newblock Mean field games with a dominating player.
\newblock {\em Applied Mathematics \& Optimization}, 74(1):91--128, 2016.

\bibitem{blum2014learning}
Avrim Blum, Nika Haghtalab, and Ariel~D Procaccia.
\newblock Learning optimal commitment to overcome insecurity.
\newblock {\em Advances in Neural Information Processing Systems}, 27, 2014.

\bibitem{campbell2021deep}
Steven Campbell, Yichao Chen, Arvind Shrivats, and Sebastian Jaimungal.
\newblock Deep learning for principal-agent mean field games.
\newblock {\em arXiv preprint arXiv:2110.01127}, 2021.

\bibitem{carmona2016finite}
Rene Carmona and Peiqi Wang.
\newblock Finite state mean field games with major and minor players.
\newblock {\em arXiv preprint arXiv:1610.05408}, 2016.

\bibitem{chavdarova2019reducing}
Tatjana Chavdarova, Gauthier Gidel, Fran{\c{c}}ois Fleuret, and Simon
  Lacoste-Julien.
\newblock Reducing noise in {GAN} training with variance reduced extragradient.
\newblock {\em Advances in Neural Information Processing Systems}, 32, 2019.

\bibitem{coniglio2020computing}
Stefano Coniglio, Nicola Gatti, and Alberto Marchesi.
\newblock Computing a pessimistic {Stackelberg} equilibrium with multiple
  followers: The mixed-pure case.
\newblock {\em Algorithmica}, 82(5):1189--1238, 2020.

\bibitem{conitzer2002complexity}
Vincent Conitzer and Tuomas Sandholm.
\newblock Complexity of mechanism design.
\newblock {\em arXiv preprint cs/0205075}, 2002.

\bibitem{conitzer2006computing}
Vincent Conitzer and Tuomas Sandholm.
\newblock Computing the optimal strategy to commit to.
\newblock In {\em Proceedings of the 7th ACM Conference on Electronic
  Commerce}, pages 82--90, 2006.

\bibitem{elie2021mean}
Romuald Elie, Emma Hubert, Thibaut Mastrolia, and Dylan Possama{\"\i}.
\newblock Mean--field moral hazard for optimal energy demand response
  management.
\newblock {\em Mathematical Finance}, 31(1):399--473, 2021.

\bibitem{elie2019tale}
Romuald Elie, Thibaut Mastrolia, and Dylan Possama{\"\i}.
\newblock A tale of a principal and many, many agents.
\newblock {\em Mathematics of Operations Research}, 44(2):440--467, 2019.

\bibitem{fiez2019convergence}
Tanner Fiez, Benjamin Chasnov, and Lillian~J Ratliff.
\newblock Convergence of learning dynamics in {Stackelberg} games.
\newblock {\em arXiv preprint arXiv:1906.01217}, 2019.

\bibitem{fu2020mean}
Guanxing Fu and Ulrich Horst.
\newblock Mean-field leader-follower games with terminal state constraint.
\newblock {\em SIAM Journal on Control and Optimization}, 58(4):2078--2113,
  2020.

\bibitem{gidel2018variational}
Gauthier Gidel, Hugo Berard, Ga{\"e}tan Vignoud, Pascal Vincent, and Simon
  Lacoste-Julien.
\newblock A variational inequality perspective on generative adversarial
  networks.
\newblock {\em arXiv preprint arXiv:1802.10551}, 2018.

\bibitem{guo2019learning}
Xin Guo, Anran Hu, Renyuan Xu, and Junzi Zhang.
\newblock Learning mean-field games.
\newblock {\em Advances in Neural Information Processing Systems}, 32, 2019.

\bibitem{guo2020general}
Xin Guo, Anran Hu, Renyuan Xu, and Junzi Zhang.
\newblock A general framework for learning mean-field games.
\newblock {\em arXiv preprint arXiv:2003.06069}, 2020.

\bibitem{guo2022mf}
Xin Guo, Anran Hu, and Junzi Zhang.
\newblock {MF-OMO}: An optimization formulation of mean-field games.
\newblock {\em arXiv preprint arXiv:2206.09608}, 2022.

\bibitem{holmstrom1987aggregation}
Bengt Holmstrom and Paul Milgrom.
\newblock Aggregation and linearity in the provision of intertemporal
  incentives.
\newblock {\em Econometrica: Journal of the Econometric Society}, pages
  303--328, 1987.

\bibitem{huang2010large}
Minyi Huang.
\newblock Large-population lqg games involving a major player: the {Nash}
  certainty equivalence principle.
\newblock {\em SIAM Journal on Control and Optimization}, 48(5):3318--3353,
  2010.

\bibitem{huang2020mean}
Minyi Huang and Xuwei Yang.
\newblock Mean field {Stackelberg} games: State feedback equilibrium.
\newblock {\em IFAC-PapersOnLine}, 53(2):2237--2242, 2020.

\bibitem{jiang2020multi}
Suhan Jiang, Xinyi Li, and Jie Wu.
\newblock Multi-leader multi-follower {S}tackelberg game in mobile blockchain
  mining.
\newblock {\em IEEE Transactions on Mobile Computing}, 2020.

\bibitem{jin2020local}
Chi Jin, Praneeth Netrapalli, and Michael Jordan.
\newblock What is local optimality in nonconvex-nonconcave minimax
  optimization?
\newblock In {\em International Conference on Machine Learning}, pages
  4880--4889. PMLR, 2020.

\bibitem{jones1999bounded}
Bryan~D Jones.
\newblock Bounded rationality.
\newblock {\em Annual Review of Political Science}, 2(1):297--321, 1999.

\bibitem{kahneman2003maps}
Daniel Kahneman.
\newblock Maps of bounded rationality: Psychology for behavioral economics.
\newblock {\em American Economic Review}, 93(5):1449--1475, 2003.

\bibitem{korobkin2003bounded}
Russell Korobkin.
\newblock Bounded rationality, standard form contracts, and unconscionability.
\newblock {\em U. Chi. L. Rev.}, 70:1203, 2003.

\bibitem{korzhyk2011stackelberg}
Dmytro Korzhyk, Zhengyu Yin, Christopher Kiekintveld, Vincent Conitzer, and
  Milind Tambe.
\newblock Stackelberg vs. {Nash} in security games: An extended investigation
  of interchangeability, equivalence, and uniqueness.
\newblock {\em Journal of Artificial Intelligence Research}, 41:297--327, 2011.

\bibitem{lee2021fast}
Sucheol Lee and Donghwan Kim.
\newblock Fast extra gradient methods for smooth structured
  nonconvex-nonconcave minimax problems.
\newblock {\em Advances in Neural Information Processing Systems},
  34:22588--22600, 2021.

\bibitem{letchford2010computing}
Joshua Letchford and Vincent Conitzer.
\newblock Computing optimal strategies to commit to in extensive-form games.
\newblock In {\em Proceedings of the 11th ACM Conference on Electronic
  Commerce}, pages 83--92, 2010.

\bibitem{letchford2009learning}
Joshua Letchford, Vincent Conitzer, and Kamesh Munagala.
\newblock Learning and approximating the optimal strategy to commit to.
\newblock In {\em International Symposium on Algorithmic Game Theory}, pages
  250--262. Springer, 2009.

\bibitem{lewis2014computational}
Richard~L Lewis, Andrew Howes, and Satinder Singh.
\newblock Computational rationality: Linking mechanism and behavior through
  bounded utility maximization.
\newblock {\em Topics in Cognitive Science}, 6(2):279--311, 2014.

\bibitem{lin2020gradient}
Tianyi Lin, Chi Jin, and Michael Jordan.
\newblock On gradient descent ascent for nonconvex-concave minimax problems.
\newblock In {\em International Conference on Machine Learning}, pages
  6083--6093. PMLR, 2020.

\bibitem{lin2020near}
Tianyi Lin, Chi Jin, and Michael~I Jordan.
\newblock Near-optimal algorithms for minimax optimization.
\newblock In {\em Conference on Learning Theory}, pages 2738--2779. PMLR, 2020.

\bibitem{luenberger1984linear}
David~G Luenberger and Yinyu Ye.
\newblock {\em Linear and Nonlinear Programming}, volume~2.
\newblock Springer, 1984.

\bibitem{mastrolia2022agency}
Thibaut Mastrolia and Jiacheng Zhang.
\newblock Agency problem and mean field system of agents with moral hazard,
  synergistic effects and accidents.
\newblock {\em arXiv preprint arXiv:2207.11087}, 2022.

\bibitem{moon2018linear}
Jun Moon and Tamer Ba{\c{s}}ar.
\newblock Linear quadratic mean field {Stackelberg} differential games.
\newblock {\em Automatica}, 97:200--213, 2018.

\bibitem{nguyen2012linear}
Son~Luu Nguyen and Minyi Huang.
\newblock Linear-quadratic-gaussian mixed games with continuum-parametrized
  minor players.
\newblock {\em SIAM Journal on Control and Optimization}, 50(5):2907--2937,
  2012.

\bibitem{nourian2013}
Mojtaba Nourian and Peter~E Caines.
\newblock $\epsilon$-{Nash} mean field game theory for nonlinear stochastic
  dynamical systems with major and minor agents.
\newblock {\em SIAM Journal on Control and Optimization}, 51(4):3302--3331,
  2013.

\bibitem{peng2019learning}
Binghui Peng, Weiran Shen, Pingzhong Tang, and Song Zuo.
\newblock Learning optimal strategies to commit to.
\newblock In {\em Proceedings of the AAAI Conference on Artificial
  Intelligence}, volume~33, pages 2149--2156, 2019.

\bibitem{perolat2021scaling}
Julien Perolat, Sarah Perrin, Romuald Elie, Mathieu Lauri{\`e}re, Georgios
  Piliouras, Matthieu Geist, Karl Tuyls, and Olivier Pietquin.
\newblock Scaling up mean field games with online mirror descent.
\newblock {\em arXiv preprint arXiv:2103.00623}, 2021.

\bibitem{perrin2020fictitious}
Sarah Perrin, Julien P{\'e}rolat, Mathieu Lauri{\`e}re, Matthieu Geist, Romuald
  Elie, and Olivier Pietquin.
\newblock Fictitious play for mean field games: Continuous time analysis and
  applications.
\newblock {\em Advances in Neural Information Processing Systems},
  33:13199--13213, 2020.

\bibitem{sannikov2008continuous}
Yuliy Sannikov.
\newblock A continuous-time version of the principal-agent problem.
\newblock {\em The Review of Economic Studies}, 75(3):957--984, 2008.

\bibitem{simon1955behavioral}
Herbert~A Simon.
\newblock A behavioral model of rational choice.
\newblock {\em The Quarterly Journal of Economics}, 69(1):99--118, 1955.

\bibitem{tambe2011security}
Milind Tambe.
\newblock {\em Security and game theory: algorithms, deployed systems, lessons
  learned}.
\newblock Cambridge University Press, 2011.

\bibitem{trejo2018adapting}
Kristal~K Trejo, Julio~B Clempner, and Alexander~S Poznyak.
\newblock Adapting attackers and defenders patrolling strategies: A
  reinforcement learning approach for {Stackelberg} security games.
\newblock {\em Journal of Computer and System Sciences}, 95:35--54, 2018.

\bibitem{yang2021linear}
Xuwei Yang and Minyi Huang.
\newblock Linear quadratic mean field {Stackelberg} games: Master equations and
  time consistent feedback strategies.
\newblock In {\em 2021 60th IEEE Conference on Decision and Control (CDC)},
  pages 171--176. IEEE, 2021.

\bibitem{zheng2020ai}
Stephan Zheng, Alexander Trott, Sunil Srinivasa, Nikhil Naik, Melvin Gruesbeck,
  David~C Parkes, and Richard Socher.
\newblock The {AI} economist: Improving equality and productivity with
  {AI}-driven tax policies.
\newblock {\em arXiv preprint arXiv:2004.13332}, 2020.

\bibitem{zhong2021can}
Han Zhong, Zhuoran Yang, Zhaoran Wang, and Michael~I Jordan.
\newblock Can reinforcement learning find {Stackelberg}-{Nash} equilibria in
  general-sum {Markov} games with myopic followers?
\newblock {\em arXiv preprint arXiv:2112.13521}, 2021.

\end{thebibliography}
\end{document}